\documentclass{amsproc}
\usepackage{amssymb}
\usepackage{amsmath}
\usepackage[matrix,arrow,curve,cmtip]{xy}

\numberwithin{equation}{subsection}

\theoremstyle{plain}
\newtheorem{bigthm}{Theorem}

\newtheorem{theorem}[equation]{Theorem}
\newtheorem{cor}[equation]{Corollary}     
\newtheorem{lemma}[equation]{Lemma}         
\newtheorem{sublemma}[equation]{Sublemma}         
\newtheorem{prop}[equation]{Proposition} 
\newtheorem{addendum}[equation]{Addendum}

\theoremstyle{definition}
\newtheorem{definition}[equation]{Definition}

\theoremstyle{remark}
\newtheorem{remark}[equation]{Remark}

\begin{document}

\title{The de~Rham-Witt complex and $p$-adic vanishing cycles${}^{^{1,2}}$}

\author{Thomas Geisser}
\address{University of Southern California, Los Angeles, California}
\email{geisser@math.usc.edu}

\author{Lars Hesselholt}
\address{Massachusetts Institute of Technology, Cambridge, Massachusetts}
\email{larsh@math.mit.edu}

\address{Nagoya University, Nagoya, Japan}
\email{larsh@math.nagoya-u.ac.jp}

\thanks{${}^{1}$ A previous version of this paper was entitled
\textsl{On the $K$-theory of a henselian discrete valuation field with
non-perfect residue field}.}

\thanks{${}^{2}$\ The authors were supported in part by grants from the
National Science Foundation. The first author received additional
support from the Japan Society for the Promotion of Science and the
Alfred P. Sloan Foundation.}

\begin{abstract}We determine the structure of the reduction modulo $p$
of the absolute de~Rham-Witt complex~\cite{hm3,hm2} of a smooth scheme
over a discrete valuation ring of mixed characteristic $(0,p)$ with
log-poles along the special fiber and show that the sub-sheaf fixed by
the Frobenius is isomorphic to the sheaf of $p$-adic vanishing
cycles. We use this result together with the main results
of~\emph{op.~cit.} to evaluate the algebraic $K$-theory with finite
coefficients of the quotient field of the henselian local ring at a
generic point of the special fiber. The result affirms the
Lichtenbaum-Quillen conjecture for this field. 
\end{abstract}



\maketitle

\section*{Introduction}

Let $V$ be a henselian discrete valuation ring with quotient field
$K$ of characteristic $0$ and perfect residue field $k$ of odd
characteristic $p$. Let $X$ be a smooth $V$-scheme of relative
dimension $r$, and let $i$ and $j$ denote the inclusion of the special
and generic fiber, respectively, as in the cartesian diagram
$$\xymatrix{
{Y\;} \ar@{^{(}->}[r]^{i} \ar[d] &
{X} \ar[d] &
{\;U} \ar@{_{(}->}[l]_{j} \ar[d] \cr
{\operatorname{Spec} k\;} \ar@{^{(}->}[r] &
{\operatorname{Spec} V} &
{\;\operatorname{Spec} K.} \ar@{_{(}->}[l] \cr
}$$
The henselian local ring of $X$ at a generic point of $Y$ is a
henselian discrete valuation ring $\mathcal{V}$ whose residue
field $\kappa$ is the non-perfect function field of a connected
component of $Y$. Let $\mathcal{K}$ be the quotient field of
$\mathcal{V}$.

We consider the ring $\mathcal{O}_X$ with the log-structure
$\alpha \colon M_X \to \mathcal{O}_X$ determined by the special
fiber. The absolute de~Rham-Witt complex of $(X,M_X)$~\cite{hm3,hm2}
is defined to be the universal Witt complex over $(\mathcal{O}_X,M_X)$
and is denoted
$$W_n \Omega_{(X,M_X)}^* = W_n \Omega_{(\mathcal{O}_X,M_X)}^*.$$
It has a natural descending filtration by the differential graded ideals
$$\operatorname{Fil}_U^m W_n\Omega_{(X,M_X)}^* \subset
W_n \Omega_{(X,M_X)}^*$$
generated by $W_n(\mathfrak{m}^j\mathcal{O}_X)$, if $m = 2j$ is even,
and by $W_n(\mathfrak{m}^j\mathcal{O}_X) \cdot d\log_n M_X$ and by
$W_n(\mathfrak{m}^{j+1}\mathcal{O}_X)$, if $m = 2j + 1$ is odd. Here
and throughout $\mathfrak{m}$ denotes the maximal ideal of $V$. There
is a natural isomorphism
$$i^* W_n\Omega_{(X,M_X)}^* /
i^*\operatorname{Fil}_U^1 W_n\Omega_{(X,M_X)}^*
\xrightarrow{\sim} W_n\Omega_Y^*$$
onto the classical de~Rham-Witt complex of
Bloch-Deligne-Illusie~\cite{illusie} of the special fiber, and
there is a natural isomorphism
$$i^* W_n\Omega_{(X,M_X)}^* /
i^*\operatorname{Fil}_U^2 W_n\Omega_{(X,M_X)}^*
\xrightarrow{\sim} W_n\Omega_{(Y,M_Y)}^*$$
onto the de~Rham-Witt complex of Hyodo-Kato~\cite{hyodokato} of the
special fiber with the induced log-structure. (The latter was denoted
by $W_n\tilde{\omega}_Y^q$ in~\emph{op.cit.}) But the full complex
$W_n\Omega_{(X,M_X)}^*$ and its quotients by the higher terms of the
filtration have not been considered before. In Sect.~\ref{drw} below
we evaluate the graded pieces for the induced filtration of the
reduced sheaves
$$E_n^q = i^*W_n \Omega_{(X,M_X)}/pi^*W_n \Omega_{(X,M_X)}^q.$$
These are quasi-coherent $E_n^0$-modules on the small \'{e}tale site
of $Y$ which we analyze as follows. A choice of local coordinates of
an open neighborhood of $X$ around a point of $Y$ determines a ring
homomorphism $\delta_n \colon \mathcal{O}_Y \to E_n^0$ defined in the
corresponding open neighborhood of $Y$. Hence, in this open
neighborhood of $Y$, we may consider the sheaves $E_n^q$ and the
graded pieces $\operatorname{gr}_U^mE_n^q$ as quasi-coherent
$\mathcal{O}_Y$-modules. We show in Thm.~\ref{Ubasis} below that both
are free $\mathcal{O}_Y$-modules and give an explicit basis. We note
here that the rank of the free $\mathcal{O}_Y$-module $E_n^q$ is
$$\operatorname{rk}_{\mathcal{O}_Y} E_n^q =
\binom{r+1}{q} \cdot e \cdot \sum_{s = 0}^{n-1} p^{rs}$$ 
and that the length of the filtration is $2e' = 2pe/(p-1)$. In
Sect.~\ref{vanishing} we use this result to analyze the subsheaf
$M_n^q$ of $E_n^q$ that is fixed by the Frobenius. The sheaf $M_n^q$
is defined to be the kernel of the map
$$R - F \colon E_n^q \to E_{n-1}^q$$
and is a sheaf of abelian groups on the small \'{e}tale site of
$Y$ in the \'{e}tale topology. The filtration of $E_n^q$ induces a
filtration of $M_n^q$ and in Thm.~\ref{gr} below we evaluate the
graded pieces $\operatorname{gr}_U^m M_n^q$ considered as sheaves of
pro-abelian groups as $n$ varies. The result agrees with the
calculation by Bloch-Kato~\cite[Cor.~1.4.1]{blochkato} of the graded
pieces $\operatorname{gr}_U^mi^*R^qj_*\mu_p^{\otimes q}$ of a
corresponding filtration of the sheaf of $p$-adic vanishing cycles. By
combining the two results we obtain the following theorem.

\begin{bigthm}\label{main}Suppose $\mu_{p^v}\subset K$. Then there
is a natural exact sequence
$$0 \to i^*R^qj_*\mu_{p^v}^{\otimes q} \to
i^*W_{\boldsymbol{\cdot}} \Omega_{(X,M_X)}^q/p^v \xrightarrow{1-F}
i^*W_{\boldsymbol{\cdot}} \Omega_{(X,M_X)}^q/p^v \to 0$$
of sheaves of pro-abelian groups on the small \'{e}tale site of $Y$ in
the \'{e}tale topology.
\end{bigthm}

We expect that Thm.~\ref{main} is valid also if $K$ does not contain 
the $p^v$th roots of unity. More precisely, we expect that the terms
in the sequence satisfy Galois descent for the extension
$K(\mu_{p^v})/K$; compare~\cite[Thm.~1(1)]{kahn}. At present, however,
the structure of the sheaves $\smash{i^*W_n\Omega_{(X,M_X)}^q/p^v}$ is
known only if $K$ contains the $p^v$th roots of unity. It appears to
be an important problem to determine the structure of these sheaves in
general.

We remark that Thm.~\ref{main} is not valid if the absolute
de~Rham-Witt complex is replaced by the relative de~Rham-Witt complex
of Langer-Zink~\cite{langerzink}. The reader is referred to~\cite{h3}
for a comparison of the two complexes.

The algebraic $K$-theory with $\mathbb{Z}/p^v$-coefficients of the
field $K$ was determined in~\cite{hm2}. In the final Sect.~\ref{bl} we
combine Thm.~\ref{main} and the main results of~\cite{hm3,hm2} to
extend this result to the field $\mathcal{K}$. Indeed we prove the
following formula predicted by the Beilinson-Lichtenbaum
conjectures~\cite{beilinson, beilinson1,lichtenbaum1}.

\begin{bigthm}\label{ktheory}Suppose that $\mu_{p^v}\subset K$. Then
the canonical map
$$K_*^M(\mathcal{K})\otimes S_{\mathbb{Z}/p^v}(\mu_{p^v}) \to
K_*(\mathcal{K},\mathbb{Z}/p^v)$$
is an isomorphism.
\end{bigthm}

The second tensor factor on the left-hand side is the symmetric
algebra on the $\mathbb{Z}/p^v$-module $\mu_{p^v}$ which is free of
rank one and the map of the statement takes a generator
$\zeta\in\mu_{p^v}$ to the associated Bott element $b_{\zeta} \in
K_2(\mathcal{K},\mathbb{Z}/p^v)$. The Milnor groups
$K_q^M(\mathcal{K})/pK_q^M(\mathcal{K})$ were evaluated by Kato 
in~\cite[Thm.~2(1)]{kato1}. They are concentrated in degrees
$0 \leqslant q\leqslant  r+2$. Hence Thm.~\ref{ktheory} shows in
particular that the algebraic $K$-groups
$K_*(\mathcal{K},\mathbb{Z}/p^v)$ are two-periodic above this range of
degrees.

The results of this paper were reported in expository form in~\cite{h2}.

In this paper, a pro-object of a category $\mathcal{C}$ will be taken
to mean a functor from the set of positive integers, viewed as a
category with one arrow from $n+1$ to $n$, to $\mathcal{C}$, and a
\emph{strict} map between pro-objects a natural transformation. A
general map between pro-objects $X$ and $Y$ of $\mathcal{C}$ is
defined to be an element of
$$\operatorname{Hom}_{\operatorname{pro}-\mathcal{C}}(X,Y) =
\operatornamewithlimits{lim}_n \operatornamewithlimits{colim}_m 
\operatorname{Hom}_{\mathcal{C}}(X_m,Y_n).$$
We view objects of $\mathcal{C}$ as constant pro-objects of
$\mathcal{C}$. Throughout, the prime $p$ is a fixed odd prime. We
abbreviate $e''=e/(p-1)$ and $e'=pe/(p-1)$.

This paper was written, in part, while the authors visited the
University of Tokyo. The second author also visited the Isaac Newton
Institute of Mathematical Sciences and Stanford University. We would
like to express our sincere gratitude for the financial support and
the hospitality that we received. We are particularly grateful to
Viorel Costeanu for allowing us include his proof of the Steinberg
relation in the de~Rham-Witt complex as an appendix to this paper.

\section{The absolute de~Rham-Witt complex}\label{drw}

\subsection{}We consider the de~Rham-Witt complex of
log-$\mathbb{Z}_{(p)}$-algebras introduced in~\cite[Sect.~3]{hm2}; see
also~\cite{hm3}. It generalizes the de~Rham-Witt complex of
log-$\mathbb{F}_p$-algebras of Hyodo-Kato~\cite{hyodokato}.

A log-ring $(A,M_A)$ in the sense of~\cite{kato} is a ring $A$ (in a
topos) with a pre-log-structure defined to be a map of monoids
$\alpha\colon M_A\to (A,\cdot\,)$ from a monoid $M_A$ to the
underlying multiplicative monoid of $A$. A log-differential graded
ring $(D,M_D)$ is a differential graded ring $D$, a pre-log-structure
$\alpha \colon M_D\to (D^0,\cdot\,)$, and a map of monoids $d\log
\colon M_D\to (D^1,+)$ that satisfies that $d\circ d\log=0$ and that
for all $x\in M_D$, $d\alpha(x)=\alpha(x)d\log x$. Maps of log-rings
and log-differential graded rings are defined in the obvious way.

Let $W_n(A)$ be the ring of Witt vectors of length $n$ in $A$. If
$\alpha \colon M_A\to A$ is a pre-log-structure, then so is the
composite
$$M_A \xrightarrow{\alpha} A \to W_n(A),$$
where the right-hand map is the multiplicative section
$a\mapsto [a]_n=(a,0,\dots,0)$. We denote this log-ring by
$(W_n(A),M_A)$. By a Witt complex over $(A,M_A)$ we mean the following
structure:

(i) a pro-log-differential graded ring
$(E_{\boldsymbol{\cdot}}^*,M_E)$ and a strict map of pro-log-rings
$$\lambda \colon (W_{\boldsymbol{\cdot}}(A),M_A) \to
(E_{\boldsymbol{\cdot}}^0,M_E);$$

(ii) a strict map of pro-log-graded rings
$$F \colon (E_{\boldsymbol{\cdot}}^*,M_E) \to
(E_{\boldsymbol{\cdot}-1}^*,M_E)$$
such that $F\lambda=\lambda F$ and such that
$$\begin{aligned}
Fd\log_n\lambda(a) & = d\log_{n-1}\lambda(a),
\hskip18.5mm\text{for all $a\in M_A$}, \cr 
Fd\lambda([a]_n) & = \lambda([a]_{n-1})^{p-1}d\lambda([a]_{n-1}),
\hskip3mm\text{for all $a\in A$;} \cr
\end{aligned}$$

(iii) a strict map of pro-graded modules over the pro-graded ring
$E^*_{\boldsymbol{\cdot}}$
$$V \colon F_*E_{\boldsymbol{\cdot}-1}^*\to E_{\boldsymbol{\cdot}}^*$$
such that $V\lambda = \lambda V$, $FV=p$, and $FdV=d$.

A map of Witt complexes over $(A,M)$ is a strict map of pro-log
differential graded rings which commutes with the maps $\lambda$, $F$
and $V$. We write $R$ for the structure map in the pro-system
$E_{\boldsymbol{\cdot}}^*$ and call it the restriction map. The
definining relations imply that $dF=pFd$ and $Vd=pdV$, but in general
there is no formula for $VF$; see~\cite[Lemma~1.2.1]{hm3}. The
de~Rham-Witt complex
$$W_{\boldsymbol{\cdot}} \Omega_{(A,M_A)}^*$$
is defined to be the initial Witt complex over $(A,M_A)$. The proof
that it exists is given in~\cite[Thm.~A]{hm3}. The proof also shows
that the canonical map
$$\lambda \colon \Omega_{(W_n(A),M_A)}^q \to
W_n \Omega_{(A,M_A)}^q$$
is surjective. Hence, every element on the right-hand side can be
written non-uniquely as a differential $q$-form on $(W_n(A),M_A)$. The
descending filtration of the de~Rham-Witt complex by the differential
graded ideals
$$\operatorname{Fil}^sW_n \Omega_{(A,M_A)}^* =
V^sW_{n-s} \Omega_{(A,M_A)}^* + dV^sW_{n-s} \Omega_{(A,M_A)}^*$$
is called the standard filtration. It satisfies that
$$\begin{aligned}
F(\operatorname{Fil}^sW_n \Omega_{(A,M_A)}^q) & \subset
\operatorname{Fil}^{s-1}W_{n-1} \Omega_{(A,M_A)}^q \cr
V(\operatorname{Fil}^sW_n \Omega_{(A,M_A)}^q) & \subset
\operatorname{Fil}^{s+1}W_{n+1} \Omega_{(A,M_A)}^q \cr
\end{aligned}$$
but, in general, it is not multiplicative. The restriction induces an
isomorphism
$$W_n \Omega_{(A,M_A)}^q/\operatorname{Fil}^sW_n \Omega_{(A,M_A)}^q
\xrightarrow{\sim}W_s \Omega_{(A,M_A)}^q.$$

\subsection{}Let $X$ be as in the introduction. We recall
from~\cite{kato} that the canonical log-structure on $X$ is given by
the cartesian square of sheaves of monoids
$$\xymatrix{
{M_X} \ar[r]^{\alpha} \ar[d] &
{\mathcal{O}_X} \ar[d] \cr
{j_*\mathcal{O}_{U}^*} \ar[r] &
{j_*\mathcal{O}_{U}.} \cr
}$$
A choice of uniformizer $\pi$ of $V$ gives rise to an isomorphism
$$\mathcal{O}_X^*\times\mathbb{N}_0 \xrightarrow{\sim} M_X$$
that takes $(u,i)$ to $\pi^iu$. In this case, the de~Rham-Witt complex
$$W_{\boldsymbol{\cdot}} \Omega_{(X,M_X)}^* =
W_{\boldsymbol{\cdot}} \Omega_{(\mathcal{O}_X,M_X)}^*$$
has an additional filtration by the differential graded ideals
$$\operatorname{Fil}_U^m W_n \Omega_{(X,M_X)}^*$$
generated by $W_n(\mathfrak{m}^j\mathcal{O}_X)$, if $m=2j$ is even,
and by $W_n(\mathfrak{m}^j\mathcal{O}_X)\cdot d\log_nM_X$ and
$W_n(\mathfrak{m}^{j+1}\mathcal{O}_X)$, if $m=2j+1$ is odd. Here
$\mathfrak{m}$ is the maximal ideal of $V$. We call this filtration
the $U$-filtration.

\begin{lemma}\label{reduction}The $U$-filtration is multiplicative and
is preserved by the restriction, Frobenius, and Verschiebung
maps. Moreover, if $X_j = X \times_{\operatorname{Spec} V}
\operatorname{Spec}(V/\mathfrak{m}^j)$, if $i_j \colon X_j\to X$ is
the closed immersion, and if $\alpha \colon M_{X_j} \to
\mathcal{O}_{X_j}$ is the induced pre-log structure, then the
canonical projection induces an isomorphism
$$i_j^*W_n \Omega_{(X,M_X)}^q / i_j^*\operatorname{Fil}_U^{2j}W_n
\Omega_{(X,M_X)}^q \xrightarrow{\sim}
W_n \Omega_{(X_j,M_{X_j})}^q.$$
\end{lemma}

\begin{proof}A functor which has a right adjoint preserves initial
objects. Hence
$$i_j^*W_n \Omega_{(X,M_X)}^q / i_j^*\operatorname{Fil}_U^{2j}W_n
\Omega_{(X,M_X)}^q = W_n \Omega_{(i_j^*\mathcal{O}_X,i_j^*M_X)}^q /
\operatorname{Fil}_U^{2j}W_n \Omega_{(i_j^*\mathcal{O}_X,i_j^*M_X)}^q.$$
Let $(B,M)$ be a log-ring, let $J\subset B$ be an ideal, and let
$(\bar{B},\bar{M})$ be the ring $\bar{B}=B/J$ with the induced
pre-log-structure given by the composition
$$\bar{\alpha} \colon \bar{M}=M \xrightarrow{\alpha}
B \xrightarrow{\operatorname{pr}} B/J.$$
One shows as in~\cite[Lemma~2.4]{gh1} that the canonical projection
induces a surjection
$$W_n \Omega_{(B,M)}^* \to W_n \Omega_{(\bar{B},\bar{M})}^*$$
and that the kernel is equal to the differential graded ideal
generated by the ideal $W_n(J)\subset W_n(B)$. The lemma is a special
case of this statement. 
\end{proof}

\begin{lemma}\label{pU}Let $e$ be the ramification index of $V$ and
$e''=e/(p-1)$. Then
$$\begin{aligned}
p \operatorname{Fil}_U^{2j}W_n \Omega_{(X,M_X)}^q {} & \subset
\operatorname{Fil}_U^{2\min\{j+e,pj\}}W_n \Omega_{(X,M_X)}^q,
\hskip 4 mm \text{for $j\geqslant 0$,} \cr
p \operatorname{Fil}_U^{2j}W_n \Omega_{(X,M_X)}^q {} & =
\operatorname{Fil}_U^{2(j+e)}W_n \Omega_{(X,M_X)}^q,
\hskip 13.4 mm \text{for $j\geqslant e''$.}\cr
\end{aligned}$$
\end{lemma}

\begin{proof}By the definition of the $U$-filtration, it suffices to
show that
$$\begin{aligned}
pW_n(\mathfrak{m}^j\mathcal{O}_X) {} & \subset
W_n(\mathfrak{m}^{\min\{j+e,pj\}}\mathcal{O}_X), 
\hskip 4 mm \text{for $j\geqslant 0$,} \cr
W_n(\mathfrak{m}^{j+e}\mathcal{O}_X) {} & \subset
pW_n(\mathfrak{m}^j\mathcal{O}_X),
\hskip 16.85 mm\text{for $j\geqslant e''$.} \cr \cr
\end{aligned}$$
Let $\pi$ be a uniformizer of $V$ with minimal polynomial
$x^e+p\theta(x)$ and recall from the proof of~\cite[Prop.~3.1.5]{hm2}
that $[\pi]^e+\theta([\pi])V(1)$ is contained in $pW_n(\mathcal{O}_X)$.
The second inclusion follows by iterated use of this
congruence. Finally, we recall from the proof
of~\cite[Lemma~3.1.1]{hm2}, that $p$ is congruent to $[p]+V(1)$ modulo
$pVW_n(\mathcal{O}_X)$. The first inclusion follows by induction,
since $p$ has valuation $e$.
\end{proof}

The map $d\log_n \colon M_X\to W_n \Omega_{(X,M_X)}^1$ gives rise to a map
of graded rings
$$T_{\mathbb{Z}}(M_X^{\operatorname{gp}}) \to W_n \Omega_{(X,M_X)}^*$$
from the tensor algebra of the group completion of the monoid
$M_X$. There is a descending filtration of the left-hand side by
graded ideals
$$\operatorname{Fil}_U^m T_{\mathbb{Z}}(M_X^{\operatorname{gp}})$$
which corresponds to the $U$-filtration on the right-hand side. To
define it, we first choose a uniformizer $\pi$ of $V$ such that we
have the isomorphism 
$$\mathcal{O}_X^* \times \mathbb{Z} \xrightarrow{\sim}
M_X^{\operatorname{gp}}$$
that takes $(u,i)$ to $\pi^iu$. We define
$\operatorname{Fil}_U^m T_{\mathbb{Z}}(M_X^{\operatorname{gp}})$ to be
$T_{\mathbb{Z}}(M_X^{\operatorname{gp}})$, if $m=0$, and to be the
graded ideal generated by $(1+\mathfrak{m}^j\mathcal{O}_X) \times
\{0\} \subset M_X^{\operatorname{gp}}$, if $m=2j$ and $j>0$, by
$(1+\mathfrak{m}\mathcal{O}_X) \times \mathbb{Z} \subset
M_X^{\operatorname{gp}}$, if $m=1$, and by
$(1+\mathfrak{m}^j\mathcal{O}_X) \times \{0\} \otimes \{1\} \times
\mathbb{Z} \subset M_X^{\operatorname{gp}} \otimes
M_X^{\operatorname{gp}}$ and $(1+\mathfrak{m}^{j+1}\mathcal{O}_X)
\times \{0\} \subset M_X^{\operatorname{gp}}$, if $m=2j+1$ and $j>1$.

\begin{lemma}\label{dlog}Let $x$ be a local section of
$\mathfrak{m}^j\mathcal{O}_X$. Then
$$d\log_n(1+x) \equiv \sum_{0\leqslant s<n}dV^s([x]_{n-s})$$
modulo $\operatorname{Fil}_U^{4j}W_n\Omega_{(X,M_X)}^1$.
\end{lemma}

\begin{proof}We first show that if $R$ is a ring and $x\in R$, then
$$[1+x]_n-[1]_n \equiv \sum_{0\leqslant s<n}V^s([x]_{n-s})$$
modulo the ideal $W_n((x^2))\subset W_n(R)$. By naturality, we may
assume that $R=\mathbb{Z}[x]$. If we write
$[1+x]_n - [1]_n = (a_0,a_1,\dots,a_{n-1})$, then the statement we
wish to show is that $a_s\equiv x$ modulo $(x^2)$, for all $0
\leqslant s < n$. The statement for $s=0$ is clear. We consider the
ghost coordinate
$$(1+x)^{p^s} - 1 = a_0^{p^s} + pa_1^{p^{s-1}} + \dots +
p^{s-1}a_{s-1}^p + p^sa_s.$$ 
The left-hand side is equivalent to $p^sx$ modulo $(x^2)$, and the
right-hand side, inductively, is equivalent to $p^sa_s$ modulo
$(x^2)$. It follows that $a_s$ is equivalent to $x$ modulo $(x^2)$ as
desired. If $x$ is a local section of $\mathfrak{m}^j\mathcal{O}_X$,
we may conclude that
$$[1+x]_n - [1]_n \equiv \sum_{0\leqslant s<n}V^s([x]_{n-s})$$
modulo $W_n(\mathfrak{m}^{2j}\mathcal{O}_X)$. Differentiating this
congruence we find that
$$d([1+x]_n) \equiv \sum_{0\leqslant s<n}dV^s([x]_{n-s})$$
modulo $\operatorname{Fil}_U^{4j}W_n\Omega_{(X,M_X)}^1$. It remains
to show that the left-hand side is congruent to $d\log_n(1+x)$ modulo
$\operatorname{Fil}_U^{4j}W_n\Omega_{(X,M_X)}^1$. By definition, we
have
$$[1+x]_n d\log_n(1+x) = d([1+x]_n),$$
and $[1+x]_n$ is a unit in
$W_n(\mathcal{O}_X/\mathfrak{m}^{2j}\mathcal{O}_X)$. Therefore, it
will suffice to show that the product $([1+x]_n - [1]_n) d([1+x]_n)$
is congruent to zero modulo
$\smash{\operatorname{Fil}_U^{4j}W_n\Omega_{(X,M_X)}^1}$. But this is
so, since the two factors lie in
$\operatorname{Fil}_U^{2j}W_n\Omega_{(X,M_X)}^*$ and since the
$U$-filtration is multiplicative.
\end{proof}

Lemma~\ref{dlog} determines the value of the map $d\log_n$ modulo
higher filtration on $\operatorname{Fil}_U^m
T_{\mathbb{Z}}(M_X^{\operatorname{gp}})$, for $m\geqslant
2$. Moreover, there is a commutative square
$$\xymatrix{
{ i^*M_X^{\operatorname{gp}} /
i^*\operatorname{Fil}_U^2M_X^{\operatorname{gp}} }
\ar[r]^(.39){d\log_n} \ar[d] &
{ i^*W_n \Omega_{(X,M_X)}^1 /
  i^*\operatorname{Fil}_U^2W_n\Omega_{(X,M_X)}^1 } \ar[d] \cr
{ M_Y^{\operatorname{gp}} } \ar[r]^(.4){d\log_n} &
{ W_n \Omega_{(Y,M_Y)}^1 } \cr
}$$
where the right-hand vertical map is an isomorphism by an argument
similar to the proof of Lemma~\ref{reduction}.

\subsection{}We recall from~\cite[Lemma~7.1.2]{hm3} that the sheaf
$$\bar{W}_n \Omega_{(X,M_X)}^q =
W_n \Omega_{(X,M_X)}^q / pW_n \Omega_{(X,M_X)}^q$$
is a quasi-coherent sheaf of $\bar{W}_n(\mathcal{O}_X)$-modules on the
small \'{e}tale site of $X$. Since this sheaf is supported on $Y$, we
may as well consider the sheaf
$$E_n^q = i^*\bar{W}_n\Omega_{(X,M_X)}^q$$
of quasi-coherent $i^*\bar{W}_n(\mathcal{O}_X)$-modules on $Y$. We
show that, Zariski locally, the sheaf $E_n^q$ has a non-canonical
structure of quasi-coherent $\mathcal{O}_Y$-module. Let $\varphi$ be
the absolute Frobenius on $Y$.

\begin{lemma}\label{modulestructure}Let $x_1,\dots,x_r$ be local
coordinates of an open neighborhood in $X$ of a point of $Y$, and let
$\bar{x}_1,\dots,\bar{x}_r$ be the corresponding local coordinates on
$Y$. Then there is, in the corresponding open neighborhood of $Y$, 
a strict map of pro-rings
$$\delta \colon \mathcal{O}_Y \to
i^*\bar{W}_{\boldsymbol{\cdot}} (\mathcal{O}_X)$$
such that $\delta_n(\bar{x}_i)=[x_i]_n$, $1\leqslant i\leqslant r$,
and such that $F\delta_n=\delta_{n-1}\varphi$.
\end{lemma}

\begin{proof}The ring homomorphism
$$f \colon W(k)[x_1,\dots,x_r] \to W(k)[x_1,\dots,x_r]$$
given by the Frobenius on $W(k)$ and by 
$f(x_i)=x_i^p$, $1\leqslant i\leqslant r$, is a lifting of the Frobenius on
$k[\bar{x}_1,\dots,\bar{x}_r]$. It determines a ring homomorphism
$$s_f \colon W(k)[x_1,\dots,x_r] \to W_n(W(k)[x_1,\dots,x_r])$$
that is characterized by $w_js_f=f^j$, $0\leqslant j<n$, and, after
reduction modulo $p$, a ring homomorphism
$$\bar{s}_f \colon k[\bar{x}_1,\dots,\bar{x}_r] \to
\bar{W}_n(W(k)[x_1,\dots,x_r]).$$
We compose this map with the ring homomorphisms
$$\bar{W}_n(W(k)[x_1,\dots,x_r]) \to
\bar{W}_n(V[x_1,\dots,x_r]) \to
\bar{W}_n(\mathcal{O}_X)$$
induced from the unique ring homomorphism $W(k) \to V$ that induces
the identity on residue fields and the chosen ring homomorphism
$g \colon V[x_1,\dots,x_r]\to\mathcal{O}_X$ to get the top horizontal
map in the following diagram.
$$\xymatrix{
{ k[\bar{x}_1,\dots,\bar{x}_r] } \ar[r]^(.55){t_g} \ar[d]^{\bar g} &
{ \bar{W}_n(\mathcal{O}_X) } \ar[d]^{p_n} \cr
{ \mathcal{O}_Y } \ar@{=}[r] \ar@{-->}[ur] &
{ \mathcal{O}_Y. } \cr
}$$
The right-hand vertical map is the composition of the restriction of
Witt vectors with kernel $V\bar{W}_n(\mathcal{O}_X)$ and the canonical
projection $\bar{\mathcal{O}}_X\to\mathcal{O}_Y$ with kernel
$\mathfrak{m}\mathcal{O}_X/p\mathcal{O}_X$. Since both ideals are
nilpotent, and since the left-hand vertical map is \'{e}tale, there
exists a unique ring homomorphism
$$\delta_n \colon \mathcal{O}_Y \to \bar{W}_n(\mathcal{O}_X)$$
making the above diagram commute. Moreover, one immediately verifies
that $R\delta_n=\delta_{n-1}$ and that $\delta_n(\bar{x}_i)=[x_i]_n$,
$1\leqslant i\leqslant r$, as stated. It remains to show that
$F\delta_n=\delta_{n-1}\varphi$, or equivalently, that the right-hand
square in the following diagram commutes.
$$\xymatrix{
{ k[\bar{x}_1, \dots, \bar{x}_r] } \ar[r]^(.65){\bar{g}}
\ar[d]^{\varphi} &
{ \mathcal{O}_Y } \ar[r]^(.4){\delta_n} \ar[d]^{\varphi} &
{ \bar{W}_n(\mathcal{O}_X) } \ar[d]^{F} \cr
{ k[\bar{x}_1, \dots, \bar{x}_r] } \ar[r]^(.65){\bar{g}} &
{ \mathcal{O}_Y } \ar[r]^(.37){\delta_{n-1}} &
{ \bar{W}_{n-1}(\mathcal{O}_X). } \cr
}$$
Since the outer square commutes, and since the left-hand square is
cocartesian, it follows that there exists a map
$\delta_{n-1}' \colon \mathcal{O}_Y\to\bar{W}_{n-1}(\mathcal{O}_X)$
that makes the right-hand square commute. To show that
$\delta_{n-1}'=\delta_{n-1}$, it will suffice to show that
$$\mathcal{O}_Y \xrightarrow{\delta_{n-1}'} \bar{W}_{n-1}(\mathcal{O}_X)
\xrightarrow{p_{n-1}} \mathcal{O}_Y$$
is the identity map. This, in turn, follows from the calculation
$$p_{n-1}F\delta_n\bar{g}=\varphi p_n\delta_n\bar{g}=\varphi\bar{g},$$
since the left-hand square of the diagram above is cocartesian.
\end{proof}

\begin{prop}\label{rank}The sheaf $E_n^q =
i^*\bar{W}_n\Omega_{(X,M_X)}^q$ has, Zariski locally on $Y$, the
structure of a free $\mathcal{O}_Y$-module of rank
$$\operatorname{rk}_{\mathcal{O}_Y}E_n^q =
\binom{r+1}{q} e \cdot \sum_{s=0}^{n-1}p^{rs}.$$
\end{prop}

\begin{proof}We consider the sheaf $E_n^q$ with the Zariski locally
defined $\mathcal{O}_Y$-module structure given by
Lemma~\ref{modulestructure}. The statement of the proposition is
unchanged by \'{e}tale extensions, so we may assume that
$X = \mathbb{A}_{V}^r$. We prove the proposition in this case by
induction on $r$. The basic case $r=0$ follows
from~\cite[Prop.~3.4.1]{hm2}. In the induction step we
use~\cite[Thm.~B]{hm3} which shows that the
$\smash{\mathcal{O}_{\mathbb{A}_{k}^r}}$-module
$\smash{i^*\bar{W}_n \Omega_{(\mathbb{A}_{V}^r,M_{\mathbb{A}_{V}^r})}^q}$
is the base-change along
$\mathbb{A}_{k}^{r-1}\hookrightarrow\mathbb{A}_{k}^r$ of the
$\mathcal{O}_{\mathbb{A}_{k}^{r-1}}$-module
$$\begin{aligned}
{} &
i^*\bar{W}_n\Omega_{(\mathbb{A}_{V}^{r-1},M_{\mathbb{A}_{V}^{r-1}})}^q
\oplus \bigoplus_{s=1}^{n-1}\bigoplus_{\,0\leqslant j<p^s\,} F_*^s
(i^*\bar{W}_{n-s}
\Omega_{(\mathbb{A}_{V}^{r-1},M_{\mathbb{A}_{V}^{r-1}})}^q) \cr
{} & \;\oplus 
i^*\bar{W}_n\Omega_{(\mathbb{A}_{V}^{r-1},M_{\mathbb{A}_{V}^{r-1}})}^{q-1}
\oplus \bigoplus_{s=1}^{n-1}\bigoplus_{\,0\leqslant j<p^s\,} F_*^s
(i^*\bar{W}_{n-s}\Omega_{(\mathbb{A}_{V}^{r-1},
M_{\mathbb{A}_{V}^{r-1}})}^{q-1}) \cr  
\end{aligned}$$
where the index $0 \leqslant j < p^s$ is required to not be divisible by
$p$. By induction and by Lemma~\ref{modulestructure}, this module has
rank
$$\begin{aligned}
{} & \binom{r}{q} \cdot e \cdot \Big( \sum_{t=0}^{n-1} p^{(r-1)t}
+ \sum_{s=1}^{n-1} (p^s-p^{s-1}) p^{(r-1)s}
\sum_{t=0}^{n-1-s} p^{(r-1)t} \Big) \cr
{} & + \binom{r}{q-1} \cdot e \cdot \Big( \sum_{t=0}^{n-1}p^{(r-1)t}
+ \sum_{s=1}^{n-1} (p^s-p^{s-1}) p^{(r-1)s}
\sum_{t=0}^{n-1-s} p^{(r-1)t} \Big), \cr
\end{aligned}$$
and since
$$\binom{r}{q-1} \cdot e + \binom{r}{q} \cdot e =
\binom{r+1}{q} \cdot e,$$
it remains to show that
$$\sum_{t=0}^{n-1} p^{(r-1)t}
+ \sum_{s=1}^{n-1} (p^s-p^{s-1}) p^{(r-1)s}
\sum_{t=0}^{n-1-s} p^{(r-1)t}
= \sum_{s=0}^{n-1} p^{rs}.$$
To see this, we rewrite the first summand on the left-hand side as
$$1 + p^{r-1}\sum_{t=0}^{n-2} p^{(r-1)t},$$
and rewrite the $s$th summand on the left-hand side as
$$p^{rs} + p^{r(s+1)-1}\sum_{t=0}^{n-(s+1)-1} p^{(r-1)t}
- p^{rs-1} \sum_{t=0}^{n-s-1} p^{(r-1)t}.$$
The statement follows.
\end{proof}

\begin{definition}\label{similar}Two local sections $\omega$ and
$\omega'$ of $E_n^q$ are \emph{similar} if there exists a unit $u
\in k^*$ such that whenever both $\omega$ and $\omega'$ are local
sections of $\operatorname{Fil}_U^mE_n^q$ then $\omega - u \cdot
\omega'$ is a local section of $\operatorname{Fil}_U^{m+1}E_n^q$. We
write $\omega \doteq \omega'$ if $\omega$ and $\omega'$ are similar.
\end{definition}

We note that if $\omega$ and $\omega'$ are similar then so are
$R(\omega)$ and $R(\omega')$,  $F(\omega)$ and $F(\omega')$,
$d(\omega)$ and $d(\omega')$, and $V(\omega)$ and $V(\omega')$. Let
$\pi$ be a uniformizer of $V$. Then by the proof
of~\cite[Prop.~3.1.5]{hm2}, the following
similarity holds in the ring $E_n^0$.
$$[\pi]_n^e \doteq V(1).$$
In the statements and proofs of Prop.~\ref{standardbasis} and
Thm.~\ref{Ubasis} below we shall abbreviate $[x] = [x]_n$ and $d\log x
= d\log_n x$.

\begin{prop}\label{standardbasis}Let $x_1,\dots,x_r$ be local
coordinates of an open neighborhood of $X$ around a point of $Y$. Then
in the corresponding open neighborhood of $Y$ the sheaf $E_n^q =
i^*\bar{W}_n\Omega_{(X,M_X)}^q$ is a free $\mathcal{O}_Y$-module with
a basis given as follows. Let $0 \leqslant s < n$, let  $0\leqslant
i_1,\dots,i_r<p^s$ and let $0 \leqslant j < e$.

{\rm (i)} If not all $i_m$, $1 \leqslant m \leqslant r$, are zero, let
$v=\min\{v_p(i_1),\dots,v_p(i_r),v_p(j-e')\}$; if $v < v_p(i_m)$, for
all $1 \leqslant m \leqslant r$, then the local sections
$$\begin{aligned}
{} & V^s([x_1]^{i_1}\dots [x_r]^{i_r}[\pi]^j
d\log x_{m_1}\dots d\log x_{m_q}), \cr
{} & dV^s([x_1]^{i_1}\dots [x_r]^{i_r}[\pi]^j
d\log x_{m_1}\dots d\log x_{m_{q-1}}), \cr
\end{aligned}$$
where $1 \leqslant m_1 < \dots < m_q \leqslant r$ (resp.~$1 \leqslant
m_1 < \dots < m_{q-1} \leqslant r$), form a basis. If $v = v_p(i_m)$,
for some $1 \leqslant m \leqslant r$, let $m$ be maximal with this
property; then the local sections
$$\begin{aligned}
{} & V^s([x_1]^{i_1}\dots [x_r]^{i_r}[\pi]^j
d\log x_{m_1}\dots d\log x_{m_q}), \cr
{} & V^s([x_1]^{i_1}\dots [x_r]^{i_r}[\pi]^j
d\log x_{m_1}\dots d\log x_{m_{q-1}}d\log\pi), \cr
{} & dV^s([x_1]^{i_1}\dots [x_r]^{i_r}[\pi]^j
d\log x_{m_1}\dots d\log x_{m_{q-1}}), \cr
{} & dV^s([x_1]^{i_1}\dots [x_r]^{i_r}[\pi]^j
d\log x_{m_1}\dots d\log x_{m_{q-2}}d\log\pi), \cr
\end{aligned}$$
where $1 \leqslant m_1 < \dots < m_q \leqslant r$ (resp. $1
\leqslant m_1 < \dots < m_{q-1} \leqslant r$, resp. $1 \leqslant m_1 <
\dots < m_{q-2} \leqslant r$), and where all $m_i \neq m$ form
a basis.

{\rm (ii)} If all $i_m$, $1 \leqslant m \leqslant r$, are zero, the local
sections
$$V^s([\pi]^jd\log x_{m_1}\dots d\log x_{m_q}),$$
where $1 \leqslant m_1 < \dots < m_q \leqslant r$, and if $s >
v_p(j-e')$, the local sections
$$dV^s([\pi]^jd\log x_{m_1}\dots d\log x_{m_{q-1}}),$$
where $1\leqslant m_1<\dots<m_{q-1}\leqslant r$, and if $s \leqslant
v_p(j-e')$, the local sections
$$V^s([\pi]^jd\log x_{m_1}\dots d\log x_{m_{q-1}} d\log\pi),$$
where $1\leqslant m_1<\dots<m_{q-1}\leqslant r$, form a basis.
\end{prop}

\begin{proof}Let $\Gamma_n^q$ be the set of local sections of $E_n^q$
listed in the statement. It is clear that the cardinality of this set
is equal to
$$\Big( \binom{r}{q} + \binom{r}{q-1} \Big) e \cdot \sum_{s = 0}^{n-1}
p^{rs} = \binom{r+1}{q} e \cdot \sum_{s = 1}^{n-1} p^{rs}.$$
Hence, by Prop.~\ref{rank}, it suffices to show that $\Gamma_n^q$
generates $E_n^q$ as an $\mathcal{O}_Y$-module. To do this, we first
show that a larger set $\Gamma_n^q{}'$ of local sections generate
$E_n^q$ as an $\mathcal{O}_Y$-module and then show that the elements
of the complement $\Gamma_n^q{}' \smallsetminus \Gamma_n^q$ can be
expressed as $\mathcal{O}_Y$-linear combinations of elements of
$\Gamma_n^q$. It will suffice to show that the elements of
$\Gamma_n^q{}' \smallsetminus \Gamma_n^q$ are similar in the sense of
Def.~\ref{similar} to $\mathcal{O}_Y$-linear combinations of the
elements of $\Gamma_n^q$. Indeed, the $U$-filtration of $E_n^q$ is
finite.

We define $\Gamma_n^q{}'$ to be the set consisting of the local
sections
$$\begin{aligned}
{} & V^s([x_1]^{i_1}\dots [x_r]^{i_r}[\pi]^j
d\log x_{m_1}\dots d\log x_{m_q}), \cr
{} & dV^s([x_1]^{i_1}\dots [x_r]^{i_r}[\pi]^j
d\log x_{m_1}\dots d\log x_{m_{q-1}}), \cr
{} & V^s([x_1]^{i_1}\dots [x_r]^{i_r}[\pi]^j
d\log x_{m_1}\dots d\log x_{m_{q-1}}d\log\pi), \cr
{} & dV^s([x_1]^{i_1}\dots [x_r]^{i_r}[\pi]^j
d\log x_{m_1}\dots d\log x_{m_{q-2}}d\log\pi), \cr
\end{aligned}$$
where $0 \leqslant s < n$, where $0 \leqslant i_1, \dots, i_r < p^s$,
where $0 \leqslant j < e$, and where $1 \leqslant m_1 < \dots < m_q
\leqslant r$  (resp.~$1 \leqslant m_1 < \dots < m_{q-1} \leqslant r$,
resp.~$1 \leqslant m_1 < \dots < m_{q-2} \leqslant r$). We also
consider the larger set $\Gamma_n^q{}''$ of local sections defined
similarly except that the indices $i_1, \dots, i_r$ are allowed to
take all non-negative integer values. We claim that the set
$\Gamma_n^q{}''$ generates $E_n^q$ as a sheaf of $k$-vector
spaces. Indeed, the case $n = 1$ is standard and the general case
follows by an induction argument based on the exact sequence
$$\operatorname{Fil}^{n-1}E_n^q \to E_n^q \xrightarrow{R}
E_{n-1}^q \to 0$$ 
and the surjective map
$$V^{n-1} + dV^{n-1} \colon E_1^q \oplus E_1^{q-1} \to
 \operatorname{Fil}^{n-1}E_n^q.$$
Here the latter map is a $k$-linear map since $k$ is a perfect
field. Next, iterated use of the relations
$$\begin{aligned}
{} V^s([x_m]^{p^s} \cdot \omega) & = [x_m] \cdot V^s(\omega) \cr
{} dV^s([x_m]^{p^s} \cdot \omega) & = [x_m] \cdot dV^s(\omega) + [x_m]
\cdot V^s(d\log x_m \cdot \omega) \cr
\end{aligned}$$
shows that the local sections in $\Gamma_n^q{}'' \smallsetminus
\Gamma_n^q{}'$ can be expressed as $\mathcal{O}_Y$-linear combinations
of the local sections in $\Gamma_n^q{}'$. Hence $\Gamma_n^q{}'$
generates $E_n^q$ as an $\mathcal{O}_Y$-module. We proceed to show
that the elements of $\Gamma_n^q{}' \smallsetminus \Gamma_n^q$ are
similar to $\mathcal{O}_Y$-linear combinations of elements of
$\Gamma_n^q$ by considering several cases.

First, let $0 < s < n$, let $0 \leqslant i_1, \dots, i_r < p^s$, and
let $0 \leqslant j < e$. We assume that not all $i_m$ are zero and
that all $i_m$ satisfy $v_p(i_m) > v_p(j - e') = v$.  Then the local
sections of the form
$$\begin{aligned}
{} & V^s([x_1]^{i_r} \dots [x_r]^{i_r} d\log x_{m_1} \dots d\log
x_{m_{q-1}} d\log \pi) \cr
{} & dV^s([x_1]^{i_r} \dots [x_r]^{i_r} d\log x_{m_1} \dots d\log
x_{m_{q-2}} d\log \pi), \cr
\end{aligned}$$
where $1 \leqslant m_1 < \dots < m_{q-1} \leqslant r$ (resp.~$1
\leqslant m_1 < \dots < m_{q-2} \leqslant r$), are elements of
$\Gamma_n^q{}' \smallsetminus \Gamma_n^q$. We will argue that the top
elements are similar to $\mathcal{O}_Y$-linear combinations of
elements of $\Gamma_n^q$; the bottom elements are treated
analogously. We note that $0 \leqslant v < s$ and define $s' = s -
v$. We also define $i_m' = p^{-v}i_m$, $1 \leqslant m \leqslant r$,
and
$$j' = p^{-v}j + e(\frac{p^{-v} - 1}{p^{-1} - 1})  
= p^{-v}(j + pe ( \frac{p^v - 1}{p - 1})).$$
Then iterated use of the similarity $[\pi]^e \doteq V(1)$ shows that
$$\begin{aligned}
{} & V^s([x_1]^{i_1} \dots [x_r]^{i_r} [\pi]^j d\log x_{m_1} \dots
d\log x_{m_{q-1}} d\log \pi) \cr
{} & \doteq V^{s'}([x_1]^{i_1'} \dots [x_r]^{i_r'} [\pi]^{j'} d\log
x_{m_1} \dots d\log x_{m_{q-1}} d\log \pi). \cr
\end{aligned}$$
The integers $i_m'$, $1 \leqslant m \leqslant r$, are all divisible by
$p$, and the integer $j'$ is not divisible by $p$. Therefore, since
$d$ is a derivation, we conclude that
$$\begin{aligned}
{} & V^{s'}([x_1]^{i_1'} \dots [x_r]^{i_r'} [\pi]^{j'} d\log
x_{m_1} \dots d\log x_{m_{q-1}} d\log \pi) \cr
{} & \doteq V^{s'}d([x_1]^{i_1'} \dots [x_r]^{i_r'} [\pi]^{j'} d\log
x_{m_1} \dots d\log x_{m_{q-1}}) \cr
\end{aligned}$$
which is equal to zero, since $s' = s - v > 0$ and $Vd = pdV$. 

We next let $0 < s < n$, let $0 \leqslant i_1, \dots, i_r < p^s$, and
let  $0 \leqslant j < e$, and assume that not all $i_m$ are zero. We
further let $1 \leqslant m \leqslant r$ be maximal with the property
that $v_p(i_m) = \min\{v_p(i_1), \dots, v_p(i_r)\}$ and assume that $v
= v_p(i_m) \leqslant v_p(j - e')$. Then the local sections
$$\begin{aligned}
{} & V^s([x_1]^{i_1}\dots [x_r]^{i_r}[\pi]^j
d\log x_{m_1}\dots d\log x_{m_q}), \cr
{} & V^s([x_1]^{i_1}\dots [x_r]^{i_r}[\pi]^j
d\log x_{m_1}\dots d\log x_{m_{q-1}}d\log\pi), \cr
{} & dV^s([x_1]^{i_1}\dots [x_r]^{i_r}[\pi]^j
d\log x_{m_1}\dots d\log x_{m_{q-1}}), \cr
{} & dV^s([x_1]^{i_1}\dots [x_r]^{i_r}[\pi]^j
d\log x_{m_1}\dots d\log x_{m_{q-2}}d\log\pi), \cr
\end{aligned}$$
where $1\leqslant m_1<\dots<m_q\leqslant r$ (resp. $1\leqslant
m_1<\dots<m_{q-1}\leqslant r$, resp. $1\leqslant
m_1<\dots<m_{q-2}\leqslant r$), and where some $m_i = m$, are elements
of $\Gamma_n^q{}' \smallsetminus \Gamma_n^q$. We again argue that the
top elements are similar to $\mathcal{O}_Y$-linear combinations of
elements of $\Gamma_n^q$; the remaining elements are treated
analogously. The integers $v$ satisfies $0 \leqslant v < s$, and we
define $s'$, $i_{\ell}'$, and $j'$ as above such that we have the
similarity
$$\begin{aligned}
{} & V^s([x_1]^{i_1} \dots [x_r]^{i_r} [\pi]^j d\log x_{m_1} \dots
d\log x_{m_q}) \cr
{} & \doteq V^{s'}([x_1]^{i_1'} \dots [x_r]^{i_r'} [\pi]^{j'} d\log
x_{m_1} \dots d\log x_{m_q}). \cr
\end{aligned}$$
We claim that this local section is similar to the local section
$$\begin{aligned}
{} & \sum_{\ell \neq m} i_{\ell}' \cdot V^s([x_1]^{i_1}
\dots [x_r]^{i_r} [\pi]^j d\log x_{\ell} d\log x_{m_1} \dots
\widehat{d\log x}_m \dots d\log x_{m_q}) \cr
{} & \hskip7mm + j' \cdot V^s([x_1]^{i_1} \dots [x_r]^{i_r} [\pi]^j
d\log \pi d\log x_{m_1} \dots \widehat{d \log x}_m \dots d\log
x_{m_q}) \cr
\end{aligned}$$
and hence similar to a $\mathcal{O}_Y$-linear combination of elements
of $\Gamma_n^q$ as desired. Here the term $d\log x_m$ is
omitted. Indeed, on the one hand, the local section
$$V^{s'}d([x_1]^{i_1'} \dots [x_r]^{i_r'} [\pi]^{j'} d\log x_{m_1}
\dots \widehat{d\log x}_m \dots d\log x_{m_q})$$
is equal to zero, since $s' > 0$ and $Vd = pdV$, and on the other
hand, this local section is equal to the sum
$$\begin{aligned}
{} & \sum_{1 \leqslant \ell \leqslant r} i_{\ell}' \cdot
V^{s'}([x_1]^{i_1'} \dots [x_r]^{i_r'} [\pi]^{j'} d\log x_{\ell} d\log
x_{m_1} \dots \widehat{d\log x}_m \dots d\log x_{m_q}) \cr
{} & \hskip7mm + j' \cdot V^{s'}([x_1]^{i_1'} \dots [x_r]^{i_r'}
[\pi]^{j'} d\log \pi d\log x_{m_1} \dots \widehat{d \log x}_m \dots
d\log x_{m_q}) \cr
\end{aligned}$$
since $d$ is a derivation. The individual summands in this sum are
similar to the corresponding summands in the sum
$$\begin{aligned}
{} & \sum_{1 \leqslant \ell \leqslant r} i_{\ell}' \cdot
V^s([x_1]^{i_1} \dots [x_r]^{i_r} [\pi]^j d\log x_{\ell} d\log x_{m_1}
\dots \widehat{d\log x}_m \dots d\log x_{m_q}) \cr
{} & \hskip7mm + j' \cdot V^s([x_1]^{i_1} \dots [x_r]^{i_r}
[\pi]^j d\log \pi d\log x_{m_1} \dots \widehat{d \log x}_m \dots
d\log x_{m_q}) \cr
\end{aligned}$$
from which the desired similarity follows since $i_m'$ is not
divisible by $p$.

Finally, let $0 \leqslant s < n$ and let $0 \leqslant j < e$. If $s
\leqslant v_p(j - e')$, the local sections
$$dV^s([\pi]^j d\log x_{m_1} \dots d\log x_{m_{q-1}}),$$
where $1 \leqslant m_1 < \dots < m_{q-1} \leqslant r$, and if $s >
v_p(j - e')$, the local sections
$$V^s([\pi]^j d\log x_{m_1} \dots d\log x_{m_{q-1}} d\log \pi),$$
where $1 \leqslant m_1 < \dots < m_{q-1} \leqslant r$, are elements of
$\Gamma_n^q{}' \smallsetminus \Gamma_n^q$. In the former case, let
$$j' = p^{-s}j + e(\frac{p^{-s} - 1}{p^{-1} - 1})  
= p^{-s}(j + pe ( \frac{p^s - 1}{p - 1})).$$
We then have the similarity
$$dV^s([\pi]^j d\log x_{m_1} \dots d\log x_{m_{q-1}}) \doteq
j' \cdot V^s([\pi]^j d\log x_{m_1} \dots d\log x_{m_{q-1}} d\log
\pi)$$
which shows that these elements are similar to $\mathcal{O}_Y$-linear
combinations of elements of $\Gamma_n^q$. In the latter case, we let
$v = v_p(j - e')$ and define $s' = s - v$ and
$$j' = p^{-v}j + e(\frac{p^{-v} - 1}{p^{-1} - 1})  
= p^{-v}(j + pe ( \frac{p^v - 1}{p - 1})).$$
Then $s' > 0$ and since $Vd = pdV$ the similarity
$$V^s([\pi]^j d\log x_{m_1} \dots d\log x_{m_{q-1}} d\log \pi)
\doteq V^{s'}d([\pi]^{j'} d\log x_{m_1} \dots d\log x_{m_{q-1}})$$
shows that these elements of $\Gamma_n^q{}' \smallsetminus \Gamma_n^q$
are similar to zero. This completes the proof of the proposition.
\end{proof}

\begin{theorem}\label{Ubasis}For $0 \leqslant j < e'$, let $v$ be the
unique integer such that
$$e \, \big(\frac{p^{-v}-1}{p^{-1}-1} \big) \leqslant j <
e \, \big(\frac{p^{-(v+1)}-1}{p^{-1}-1} \big). $$
Let $x_1,\dots,x_r$ be local coordinates of an open
neighborhood of $X$ of a point of $Y$. Then, in the corresponding open
neighborhood of $Y$, the sheaf
$$E_n^q = i^*W_n\Omega_{(X,M_X)}^q / pi^*W_n\Omega_{(X,M_X)}^q$$
has the following structure. 

{\rm(i)} If $0 \leqslant v < n$, then $\operatorname{gr}_U^{2j}E_n^q$
is a free $\mathcal{O}_Y$-module with a basis given as follows. If $p$
does not divide $j$ (resp.~if $p$ divides $j$), let $0 < s < n-v$, and
let $0 \leqslant i_1,\dots,i_r < p^s$ (resp.~let $0 \leqslant
i_1,\dots,i_r < p^s$ not all divisible by $p$); let $1 \leqslant m
\leqslant r$ be maximal with $i_m$ not divisible by $p$. Then the
local sections
$$\begin{aligned}
{} & V^s([x_1]^{i_1} \dots [x_r]^{i_r} [\pi]^j d\log x_{m_1} \dots
d\log x_{m_q}) \cr
{} & dV^s([x_1]^{i_1} \dots [x_r]^{i_r} [\pi]^j d\log x_{m_1'} \dots
d\log x_{m_{q-1}'}), \cr
\end{aligned}$$
where $1 \leqslant m_1 < \dots < m_q \leqslant r$ and $1 \leqslant
m_1' < \dots < m_{q-1}' \leqslant r$, and where all $m_i$ and $m_i'$
are different from $m$, together with the local sections
$$[\pi]^j d\log x_{m_1} \dots d\log x_{m_q},$$
where $1 \leqslant m_1 < \dots < m_q \leqslant r$, form a basis.

{\rm(ii)} If $0 \leqslant v < n$, then
$\operatorname{gr}_U^{2j+1}E_n^q$ is a free $\mathcal{O}_Y$-module
with a basis given as follows. Let $0 < s < n-v$, and let $0 \leqslant
i_1,\dots,i_r < p^s$ not all divisible by $p$; let $1 \leqslant m
\leqslant r$ be maximal with $i_m$ not divisible by $p$. Then the
local sections
$$\begin{aligned}
{} & V^s([x_1]^{i_1} \dots [x_r]^{i_r} [\pi]^j d\log x_{m_1} \dots
d\log x_{m_{q-1}}d\log \pi) \cr
{} & dV^s([x_1]^{i_1} \dots [x_r]^{i_r} [\pi]^j d\log x_{m_1'} \dots
d\log x_{m_{q-2}'}d\log \pi), \cr
\end{aligned}$$
where $1 \leqslant m_1 < \dots < m_{q-1} \leqslant r$ and
$1 \leqslant m_1' < \dots < m_{q-2}' \leqslant r$, and where all $m_i$
and $m_i'$ are different from $m$, together with the local sections
$$[\pi]^j d\log x_{m_1} \dots d\log x_{m_{q-1}}d\log \pi,$$
where $1 \leqslant m_1 < \dots < m_{q-1} \leqslant r$, form a basis.

{\rm(iii)} If $n \leqslant v$ or if $e' \leqslant j$, then
$\operatorname{Fil}_U^{2j}E_n^q = 0$.
\end{theorem}

\begin{proof}Let $\Sigma_n^q$ and $\Gamma_n^q$ denote the sets of
local sections of $E_n^q$ listed in the statements of the
Thm.~\ref{Ubasis} and of Prop.~\ref{standardbasis}, respectively.
We first construct a bijection $f \colon \Gamma_n^q \xrightarrow{\sim}
\Sigma_n^q$ with the property that $\omega$ and $f(\omega)$ are
similar in the sense of Def.~\ref{similar}. This proves that
$\Sigma_n^q$ is an $\mathcal{O}_Y$-module basis of $E_n^q$.

Let $0 \leqslant s < n$, let $0 \leqslant i_1, \dots, i_r < p^s$, and
let $0 \leqslant j < e$. We let
$$v = \min\{s, v_p(i_1), \dots, v_p(i_r), v_p(j - e')\}$$
and define $s' = s - v$, $i_m' = p^{-v}i_m$, and
$$j'= p^{-v}j + e(\frac{p^{-v} - 1}{p^{-1} - 1}) =
p^{-v}(j + pe(\frac{p^v - 1}{p - 1})) =
p^{-v}(j - e') + e'.$$
Then $0 \leqslant s' < n - v$, $0 \leqslant i_m' <
p^{s'}$, and $j'$ is an integer that satisfies
$$e(\frac{p^{-v}-1}{p^{-1}-1}) \leqslant j ' < 
e(\frac{p^{-(v+1)}-1}{p^{-1}-1}).$$
The number $j'$ is an integer since $j' \in \mathbb{Z}[\frac{1}{p}]$ and
$v_p(j') \geqslant 0$.

Conversely, let $0 \leqslant j' < e'$, and let $v$ be the unique integer
given by the above pair of inequalities. Let also $0 \leqslant s' < n -
v$ and $0 \leqslant i_1', \dots, i_r' < p^{s'}$ be given. We define $s
= s' + v$, $i_m = p^v i_m'$, and
$$j = p^v(j' - e(\frac{p^{-v} - 1}{p^{-1} - 1})) =
p^vj' - pe(\frac{p^v - 1}{p - 1}) =
p^v(j' - e') + e'.$$
Then $v \leqslant s < n$, $0 \leqslant i_m < p^s$, and $j$ is an
integer and $0 \leqslant j < e$. We define the bijection $f \colon
\Gamma_n^q \to \Sigma_n^q$ by considering several cases. In each case
the similarity of $\omega$ and $f(\omega)$ follows by iterated use of
the similarity $[\pi]^e \doteq V(1)$. 

First, suppose that not all $i_m$, $1 \leqslant m \leqslant r$, are
equal zero and that $v < v_p(i_m)$, for all $1 \leqslant m \leqslant
r$. Then $f$ takes the local sections
$$\begin{aligned}
{} & V^s([x_1]^{i_1}\dots [x_r]^{i_r}[\pi]^j
d\log x_{m_1}\dots d\log x_{m_q}),\cr
{} & dV^s([x_1]^{i_1}\dots [x_r]^{i_r}[\pi]^j
d\log x_{m_1}\dots d\log x_{m_{q-1}}), \cr
\end{aligned}$$
where $1\leqslant m_1<\dots<m_q\leqslant r$ {\rm (}resp.~$1\leqslant
m_1<\dots<m_{q-1}\leqslant r$\,{\rm )}, to the similar local sections
$$\begin{aligned}
{} & V^{s'}([x_1]^{i_1'}\dots [x_r]^{i_r'}[\pi]^{j'}
d\log x_{m_1}\dots d\log x_{m_q}),\cr
{} & dV^{s'}([x_1]^{i_1'}\dots [x_r]^{i_r'}[\pi]^{j'}
d\log x_{m_1}\dots d\log x_{m_{q-1}}). \cr
\end{aligned}$$
We note that $0 < s' < n - v$, that all $i_m'$, $1 \leqslant m
\leqslant r$, are divisible by $p$ and that $j'$ is not divisible by
$p$.

Next, suppose that not all $i_m$, $1 \leqslant m \leqslant r$, are
zero and that $v = v_p(i_m)$, for some $1 \leqslant m \leqslant
r$. Let $m$ be maximal with this property. Then $f$ takes the local
sections
$$\begin{aligned}
{} & V^s([x_1]^{i_1}\dots [x_r]^{i_r}[\pi]^j
d\log x_{m_1}\dots d\log x_{m_q}), \cr
{} & V^s([x_1]^{i_1}\dots [x_r]^{i_r}[\pi]^j
d\log x_{m_1}\dots d\log x_{m_{q-1}}d\log\pi), \cr
{} & dV^s([x_1]^{i_1}\dots [x_r]^{i_r}[\pi]^j
d\log x_{m_1}\dots d\log x_{m_{q-1}}), \cr
{} & dV^s([x_1]^{i_1}\dots [x_r]^{i_r}[\pi]^j
d\log x_{m_1}\dots d\log x_{m_{q-2}}d\log\pi), \cr
\end{aligned}$$
where $1\leqslant m_1<\dots<m_q\leqslant r$ (resp. $1\leqslant 
m_1<\dots<m_{q-1}\leqslant r$, resp. $1\leqslant
m_1<\dots<m_{q-2}\leqslant r$), and where all $m_i\neq m$, to the
similar local sections
$$\begin{aligned}
{} & V^{s'}([x_1]^{i_1'}\dots [x_r]^{i_r'}[\pi]^{j'}
d\log x_{m_1}\dots d\log x_{m_q}), \cr
{} & V^{s'}([x_1]^{i_1'}\dots [x_r]^{i_r'}[\pi]^{j'}
d\log x_{m_1}\dots d\log x_{m_{q-1}}d\log\pi), \cr
{} & dV^{s'}([x_1]^{i_1'}\dots [x_r]^{i_r'}[\pi]^{j'}
d\log x_{m_1}\dots d\log x_{m_{q-1}}), \cr
{} & dV^{s'}([x_1]^{i_1'}\dots [x_r]^{i_r'}[\pi]^{j'}
d\log x_{m_1}\dots d\log x_{m_{q-2}}d\log\pi). \cr
\end{aligned}$$
We note that $0 < s' < n - v$, that $i_m'$ is not divisible by $p$,
and that $j'$ may or may not be divisible by $p$.

Next, if all $i_m$, $1 \leqslant m \leqslant r$, are zero and if $v <
s$, then $f$ takes the local sections
$$\begin{aligned}
{} & V^s([\pi]^j d\log x_{m_1} \dots d\log x_{m_q}), \cr
{} & dV^s([\pi]^j d\log x_{m_1} \dots d\log x_{m_{q-1}}), \cr
\end{aligned}$$
where $1 \leqslant m_1 < \dots < m_q \leqslant r$ (resp.~$1 \leqslant
m_1 < \dots < m_{q-1} \leqslant r$) to the similar local sections
$$\begin{aligned}
{} & V^{s'}([\pi]^{j'} d\log x_{m_1} \dots d\log x_{m_q}), \cr
{} & dV^{s'}([\pi]^{j'} d\log x_{m_1} \dots d\log x_{m_{q-1}}). \cr
\end{aligned}$$
We note that $0 < s' < n - v$ and that $j'$ is not divisible by $p$.

Finally, if all $i_m$, $1 \leqslant m \leqslant r$, are zero and if $v
= s$, then $f$ takes the local sections
$$\begin{aligned}
{} & V^s([\pi]^j d\log x_{m_1} \dots d\log x_{m_q}), \cr
{} & V^s([\pi]^j d\log x_{m_1} \dots d\log x_{m_{q-1}} d\log \pi), \cr
\end{aligned}$$
where $1 \leqslant m_1 < \dots < m_q \leqslant r$ (resp.~$1 \leqslant
m_1 < \dots < m_{q-1} \leqslant r$) to the similar local sections
$$\begin{aligned}
{} & [\pi]^{j'} d\log x_{m_1} \dots d\log x_{m_q}, \cr
{} & [\pi]^{j'} d\log x_{m_1} \dots d\log x_{m_{q-1}} d\log \pi. \cr
\end{aligned}$$
The integer $j'$ may or may not be divisible by $p$. This completes
the definition of the map $f \colon \Gamma_n^q \to \Sigma_n^q$. It is
clear that $f$ is a bijection. 

It remains to show that the $U$-filtration of $E_n^q$ is as
stated. Let $A^{m,q}\subset E_n^q$ be the sub-$\mathcal{O}_Y$-module 
generated by those elements of $\Sigma_n^q$ that are listed in the
statement of the theorem as having $U$-filtration greater than or
equal to $m$. It is then clear that $A^{m,q}\subset
U^{m,q}=\operatorname{Fil}_U^mE_n^q$ and we must show that also
$U^{m,q}\subset A^{m,q}$. We recall that if $m=2j$ (resp.~if
$m=2j+1$), then $U^{m,*}$ is the differential graded ideal generated
by $i^*\bar{W}_n(\mathfrak{m}^j\mathcal{O}_X) \subset
i^*\bar{W}_n(\mathcal{O}_X)$ (resp.~by
$i^*\bar{W}_n(\mathfrak{m}^j\mathcal{O}_X)\cdot d\log M_X \subset
i^*\bar{W}_n \Omega_{(X,M_X)}^1$ and
$i^*\bar{W}_n(\mathfrak{m}^{j+1}\mathcal{O}_X)\subset
i^*\bar{W}_n(\mathcal{O}_X)$). So it suffices to show that
$i^*\bar{W}_n(\mathfrak{m}^j\mathcal{O}_X) \subset A^{2j,0}$, that
the product takes $A^{m,q}\otimes A^{m',q'}$ to $A^{m+m',q+q'}$, and
that the differential takes $A^{m,q}$ to $A^{m,q+1}$. The second
statement is verified by explicitly calculating the products of basis
elements of $A^{m,q}$ and $A^{m',q'}$ in a manner similar
to~\cite[Sect.~4]{hm3}, and the last statement is immediate. We verify
the first statement.

One sees as in the proof of Prop.~\ref{standardbasis}
that $U^{2j,0} = i^*\bar{W}_n(\mathfrak{m}^j\mathcal{O}_X)$ is
generated as an $\mathcal{O}_Y$-module by the local sections
$$V^s([x_1]^{i_1} \dots [x_r]^{i_r} [\pi]^a)$$
where $0 \leqslant s < n$, $0 \leqslant i_1, \dots, i_r < p^s$, and $j
\leqslant a < j + e$. Iterated use of the similarity $[\pi]^e \doteq
V(1)$ shows that the generators with $e' \leqslant a < j + e$ are
equal to zero. For the remaining generators we again let
$$v = \min\{s, v_p(i_1), \dots, v_p(i_r), v_p(a - e')\}$$
and define $s' = s - v$, $i_m' = p^{-v}i_m$, and $a' = p^{-v}(a - e')
+ e'$. Then
$$V^s([x_1]^{i_1} \dots [x_r]^{i_r}[\pi]^a) \doteq
V^{s'}([x_1]^{i_1'} \dots [x_r]^{i_r'} [\pi]^{a'})$$
and hence the local sections on the right-hand side generate
$U^{2j,0}$ as an $\mathcal{O}_Y$-module. Here $0 \leqslant s' < n -
v$, $0  \leqslant i_1',\dots,i_r' < p^{s'}$, and $j \leqslant a' <
e'$, and if $s' > 0$ then not all of $i_1', \dots, i_r'$ and $a'$ are
divisible by $p$. These local sections are all contained in $A^{2j,0}$
and hence $U^{2j,0} \subset A^{2j,0}$ as desired.  
\end{proof}

\begin{addendum}\label{drwlocalization}There is a natural exact
sequence
$$0 \to i^*\bar{W}_n \Omega_X^q \xrightarrow{j_*}
i^*\bar{W}_n \Omega_{(X,M_X)}^q \xrightarrow{\partial}
\bar{W}_n \Omega_Y^{q-1} \to 0.$$
\end{addendum}

\begin{proof}The map $j_*$ is induced by the canonical map from $X$
with the trivial log-structure to $X$ with the canonical
log-structure. To construct the map $\partial$ we first show that the
map
$$f_{\pi} \colon W_n \Omega_Y^q \oplus W_n \Omega_Y^{q-1} \to
W_n \Omega_{(Y,M_Y)}^q$$
that to $(\omega,\omega')$ assigns $\omega+\omega'd\log_n\pi$ is an
isomorphism. Since the statement is local for the \'{e}tale topology,
it suffices to consider $Y=\mathbb{A}_{k}^r$. This case follows
inductively from the trivial case $r = 0$ since the domain and target
for $Y=\mathbb{A}_{k}^r$ is given by the same
formula~\cite[Thm.~B]{hm3} in terms of the domain and target for
$Y=\mathbb{A}_{k}^{r-1}$. We now define $\partial$ to be the composite
map
$$i^*W_n \Omega_{(X,M_X)}^q \to
W_n \Omega_{(Y,M_Y)}^q \xleftarrow[\sim]{f_{\pi}}
W_n \Omega_Y^q\oplus W_n \Omega_Y^{q-1} \to
W_n \Omega_Y^{q-1}$$
where the left-hand map is the canonical projection and where the
right-hand map is the projection onto the second summand. We note that
$\partial$ is independent of the choice of uniformizer since
$f_{\pi}(\omega,\omega') = f_{u\pi}(\omega + \omega'd\log_n
u,\omega')$. It also is clear from the definition that composite
$\partial \circ j_*$ is equal to zero. 

Let $x_1,\dots,x_r$ be local coordinates of an open
neighborhood on $X$ around a point of $Y$, and let
$\bar{x}_1,\dots,\bar{x}_r$ be the corresponding local coordinates
on $Y$. Lemma~\ref{modulestructure} allows us to view the sequence of
the statement as a sequence of $\mathcal{O}_Y$-modules, and
Thm.~\ref{Ubasis} gives a basis of the middle term. The map
$\partial$ takes
$$\begin{aligned}
\partial( & V^s([x_1]^{i_1}\dots [x_r]^{i_r}d\log x_{m_1}\dots
d\log x_{m_{q-1}}d\log\pi) ) \cr
{} & = V^s([\bar{x}_1]^{i_1}\dots
[\bar{x}_r]^{i_r}d\log\bar{x}_{m_1}\dots d\log\bar{x}_{m_{q-1}}), \cr
\partial( & dV^s([x_1]^{i_1}\dots [x_r]^{i_r}d\log x_{m_1}\dots
d\log_{x_{m_{q-2}}}d\log\pi) ) \cr
{} & = dV^s([\bar{x}_1]^{i_1}\dots
[\bar{x}_r]^{i_r}d\log\bar{x}_{m_1}\dots d\log\bar{x}_{m_{q-2}}), \cr
\end{aligned}$$
and annihilates all remaining basis elements. One shows, in a manner
similar to the proof of Prop.~\ref{rank} that, Zariski locally
on $Y$, the sheaves $i^*\bar{W}_n \Omega_{X}^q$ and
$\bar{W}_n \Omega_Y^{q-1}$ with the $\mathcal{O}_Y$-module structure
given by Lemma~\ref{modulestructure} are free and that their ranks
satisfy the equation
$$\operatorname{rk}_{\mathcal{O}_Y}i^*\bar{W}_n \Omega_{X}^q +
\operatorname{rk}_{\mathcal{O}_Y}\bar{W}_n \Omega_{Y}^{q-1}
= \operatorname{rk}_{\mathcal{O}_Y}i^*\bar{W}_n \Omega_{(X,M_X)}^q.$$
This completes the proof.
\end{proof}

\subsection{}We end this section with the following result on the
structure of the higher torsion in the de~Rham-Witt complex. The proof
we give here uses the cyclotomic trace; see~\cite{h2}. It would be
desirable to have a purely algebraic proof.

\begin{prop}\label{highertorsion}If $\mu_{p^v}\subset K$, then for
all $0\leqslant m<v$ and all $q\geqslant 0$, multiplication by $p^m$
induces an isomorphism of sheaves of pro-abelian groups
$$\bar{W}_{\boldsymbol{\cdot}}\Omega_{(X,M_X)}^q = 
\operatorname{gr}_p^0W_{\boldsymbol{\cdot}}\Omega_{(X,M_X)}^q
\xrightarrow{\sim}
\operatorname{gr}_p^mW_{\boldsymbol{\cdot}}\Omega_{(X,M_X)}^q.$$
\end{prop}

\begin{proof}We must show that for all $0 \leqslant m < v$ and all
$q \geqslant 0$, the following sequence of sheaves of pro-abelian
groups on the small \'{e}tale site of $X$ is exact.
$$0 \to W_{\boldsymbol{\cdot}}\Omega_{(X,M_X)}^q/p \xrightarrow{p^m}
W_{\boldsymbol{\cdot}}\Omega_{(X,M_X)}^q/p^{m+1}
\xrightarrow{\operatorname{pr}}
W_{\boldsymbol{\cdot}}\Omega_{(X,M_X)}^q/p^m \to 0.$$
This is equivalent to the statement that for all $0 \leqslant m < v$
and all $q , s \geqslant 0$, the following sequence of pro-abelian
groups is exact
$$0 \to W_{\boldsymbol{\cdot}}\Omega_{(X,M_X)}^{q-2s} \otimes
\mu_p^{\otimes s} \to
W_{\boldsymbol{\cdot}}\Omega_{(X,M_X)}^{q-2s} \otimes
\mu_{p^{m+1}}^{\otimes s} \to
W_{\boldsymbol{\cdot}}\Omega_{(X,M_X)}^{q-2s} \otimes
\mu_{p^m}^{\otimes s} \to 0.$$
We need only show that the left-hand map is a monomorphism of
pro-abelian groups. To this end we recall that for all $0 \leqslant m
< v$ and all $q \geqslant 0$,~\cite[Thm.~E]{hm3} gives an isomorphism
of pro-abelian groups
$$\bigoplus_{s \geqslant 0}
W_{\boldsymbol{\cdot}}\Omega_{(X,M_X)}^{q-2s} \otimes
\mu_{p^m}^{\otimes s} \xrightarrow{\sim}
\operatorname{TR}_q^{\,\boldsymbol{\cdot}}(X|X_K;p,\mathbb{Z}/p^m).$$
In particular, for all $0 \leqslant m < v$ and $q \geqslant 0$, the
map induced from the reduction
$$\operatorname{TR}_q^{\,\boldsymbol{\cdot}}(X|X_K;p,\mathbb{Z}/p^{m+1})
\to
\operatorname{TR}_q^{\,\boldsymbol{\cdot}}(X|X_K;p,\mathbb{Z}/p^m)$$
is an epimorphism of pro-abelian groups. It follows that the
long-exact coefficient sequence breaks up into short-exact sequences
of pro-abelian groups
$$0\to
\operatorname{TR}_q^{\,\boldsymbol{\cdot}}(X|X_K;p,\mathbb{Z}/p) \to
\operatorname{TR}_q^{\,\boldsymbol{\cdot}}(X|X_K;p,\mathbb{Z}/p^{m+1}) \to
 \operatorname{TR}_q^{\,\boldsymbol{\cdot}}(X|X_K;p,\mathbb{Z}/p^m) \to 0.$$
The proposition follows.
\end{proof}

\section{$p$-adic vanishing cycles}\label{vanishing}

\subsection{}Let $(i^*\bar{W}_n\Omega_{(X,M_X)}^q)^{F=1}$
and $(i^*\bar{W}_n\Omega_{(X,M_X)}^q)_{F=1}$ denote the kernel and
cokernel, respectively, of the map
$$R - F \colon i^*\bar{W}_n\Omega_{(X,M_X)}^q \to 
i^*\bar{W}_{n-1}\Omega_{(X,M_X)}^q$$
of sheaves of abelian groups on the small \'{e}tale site of $Y$. We
consider these sheaves both in the Nisnevich topology and the
\'{e}tale topology. The $U$-filtration is preserved by $R - F$ and
hence induces filtrations of the kernel and cokernel sheaves. We begin
with the following observation.

\begin{lemma}\label{R-F}Suppose that $m\geqslant 2$. Then, for all
integers $n$ and $q$, the map 
$$R-F \colon  \operatorname{Fil}_U^m\bar{W}_n\Omega_{(X,M_X)}^q \to
\operatorname{Fil}_U^m\bar{W}_{n-1}\Omega_{(X,M_X)}^q$$
is a surjective map of presheaves of abelian groups on the small
\'{e}tale site of $Y$.
\end{lemma}

\begin{proof}We consider the case $m=2j$; the case $m=2j+1$ is
similar. It suffices to show that if $a_0,\dots,a_q$ and
$a_0',\dots,a_{q-1}'$ are local sections of $\mathcal{O}_X$ such that
$\operatorname{ord}_Y(a_i) \geqslant j$ and
$\operatorname{ord}_Y(a_{i'}')\geqslant j$, for some $0 \leqslant
i\leqslant q$ and $0 \leqslant i' \leqslant q-1$, then the following
local sections are in the image of $R-F$.
$$\begin{aligned}
{} & V^{s_0}[a_0]_{n-1}dV^{s_1}[a_1]_{n-1}\dots dV^{s_q}[a_q]_{n-1}, \cr
{} & V^{s_0}[a_0]_{n-1}dV^{s_1}[a_1]_{n-1}\dots
dV^{s_{q-1}}[a_{q-1}]_{n-1}d\log_{n-1}\pi. \cr
\end{aligned}$$
Indeed, every local section of
$\operatorname{Fil}_U^m\bar{W}_{n-1}\Omega_{(X,M_X)}^q$ is a sum
of such elements. We now use that since $j \geqslant 1$, the
following series converge.
$$\begin{aligned}
{} & \sum_{t \geqslant 0} F^t(V^{s_0}[a_0]_{n+t}dV^{s_1}[a_1]_{n+t}
\dots dV^{s_q}[a_q]_{n+t}), \cr
{} & \sum_{t \geqslant 0} F^t(V^{s_0}[a_0]_{n+t}dV^{s_1}[a_1]_{n+t}
\dots dV^{s_{q-1}}[a_{q-1}]_{n+t}d\log_{n+t}\pi). \cr
\end{aligned}$$
The images by $R-F$ of the sums of these series are equal to the given
elements.
\end{proof}

Let $\Omega_{Y,\log}^q \subset \Omega_Y^q$ be the subsheaf generated
locally for the \'{e}tale topology on $Y$ by the local sections of the
form $d\log y_1 \dots d\log y_q$. If $y$ is a local section of
$\mathcal{O}_Y$, we denote by $\tilde{y}$ any lifting of $y$ to a
local section of $i^*\mathcal{O}_X$.

\begin{theorem}\label{gr}The sheaf $M_{\boldsymbol{\cdot}}^q =
(i^*\bar{W}_{\boldsymbol{\cdot}}\Omega_{(X,M_X)}^q)^{F=1}$ of
pro-abelian groups on the small \'{e}tale site of $Y$ has the
following structure.

{\rm(i)} There is an isomorphism
$$\rho_0 \colon \Omega_{Y,\log}^q \xrightarrow{\sim}
\operatorname{gr}_U^0M_{\boldsymbol{\cdot}}^q \hskip6mm (\text{resp.}\
\rho_1 \colon \Omega_{Y,\log}^{q-1} \xrightarrow{\sim}
\operatorname{gr}_U^1M_{\boldsymbol{\cdot}}^q)$$
that to $d\log y_1 \dots d\log y_q$ (resp.~$d\log y_1 \dots d\log
y_{q-1}$) assigns $d\log\tilde{y}_1\dots d\log\tilde{y}_q$
(resp. $d\log\tilde{y}_1 \dots d\log\tilde{y}_{q-1}d\log\pi$).

{\rm(ii)} If $0 < j < e'$, and if $p$ does not divide $j$ (resp.~if
$p$ divides $j$), there is an isomorphism
$$\rho_{2j} \colon \Omega_Y^{q-1}/B\Omega_Y^{q-1} \xrightarrow{\sim}
\operatorname{gr}_U^{2j}M_{\boldsymbol{\cdot}}^q \hskip6mm (\text{resp.}\
\rho_{2j} \colon \Omega_Y^{q-1}/Z\Omega_Y^{q-1} \xrightarrow{\sim}
\operatorname{gr}_U^{2j}M_{\boldsymbol{\cdot}}^q)$$
that to $ad\log y_1 \dots d\log y_{q-1}$ assigns
$d\log(1+\pi^j\tilde{a})d\log\tilde{y}_1 \dots d\log\tilde{y}_{q-1}$.

{\rm(iii)} If $0 < j < e'$, there is an isomorphism
$$\rho_{2j+1} \colon \Omega_Y^{q-2}/Z\Omega_Y^{q-2} \xrightarrow{\sim}
\operatorname{gr}_U^{2j+1}M_{\boldsymbol{\cdot}}^q$$
that takes $ad\log y_1 \dots d\log y_{q-2}$ to
$d\log(1+\pi^j\tilde{a})d\log\tilde{y}_1 \dots d\log\tilde{y}_{q-2}d\log\pi$.

{\rm(iv)} If $e' \leqslant j$, then
$\operatorname{Fil}_U^{2j}M_{\boldsymbol{\cdot}}^q$ is equal to zero.
\end{theorem}

\begin{proof}It follows from Prop.~\ref{steinberg} and
from~\cite[Sect.~4]{blochkato} that the maps $\rho_m$ of the statement
are well-defined strict maps of sheaves of pro-abelian groups. We first
consider the statement~(i). We abbreviate $E_n^q =
i^*\bar{W}_n\Omega_{(X,M_X)}^q$ as before. It follows from
Lemma~\ref{R-F} that there is an exact sequence of sheaves of abelian
groups
$$0 \to M_n^q/\operatorname{Fil}_U^2M_n^q \to
E_n^q/\operatorname{Fil}_U^2E_n^q \xrightarrow{R-F}
E_{n-1}^q/\operatorname{Fil}_U^2E_{n-1}^q.$$
Moreover, Lemma~\ref{reduction} identifies the middle and right-hand
terms with the reduction modulo $p$ of the de~Rham-Witt complex of
$(Y,M_Y)$. Hence, we have an isomorphism of sheaves of pro-abelian
groups
$$M_{\boldsymbol{\cdot}}^q/\operatorname{Fil}_U^2M_{\boldsymbol{\cdot}}^q
\xrightarrow{\sim} (\bar{W}_{\boldsymbol{\cdot}}\Omega_{(Y,M_Y)}^q)^{F=1}.$$
The structure of the right-hand side is well-understood;
see~\cite[Prop.~2.4.1]{tsuji}. The statement for
$\operatorname{gr}_U^0M_{\boldsymbol{\cdot}}^q$ and
$\operatorname{gr}_U^1M_{\boldsymbol{\cdot}}^q$ follows.

We next prove the statement~(ii). It follows again from
Lemma~\ref{R-F} that there is a short-exact sequence of sheaves of
abelian groups
$$0 \to \operatorname{gr}_U^{2j}M_n^q \to
\operatorname{gr}_U^{2j}E_n^q \xrightarrow{R-F}
\operatorname{gr}_U^{2j}E_{n-1}^q \to 0.$$
We recall that if $x_1,\dots,x_r$ are local coordinates of an open
neighborhood of $X$ of a point of $Y$, then, in the corresponding open
neighborhood of $Y$, the sheaf $\operatorname{gr}_U^{2j}E_n^q$ has the
structure of a free $\mathcal{O}_Y$-module with a basis given by the
local sections
$$\begin{aligned}
{} & V^s([x_1]^{i_1} \dots [x_r]^{i_r} [\pi]^j d\log x_{m_1} \dots
d\log x_{m_q}), \cr
{} & dV^s([x_1]^{i_1} \dots [x_r]^{i_r} [\pi]^j d\log x_{m_1} \dots
d\log x_{m_{q-1}}), \cr
\end{aligned}$$
where $0 \leqslant s < n - v$ and $0 < s < n - v$, respectively, and
where the multi-indices $i$ and $m$ vary as in the statement of
Thm.~\ref{Ubasis}~(ii). We consider the short-exact sequence of
$\mathcal{O}_Y$-modules
$$0 \to (\operatorname{gr}_U^{2j}E_n^q)' \to \operatorname{gr}_U^{2j}E_n^q \to (\operatorname{gr}_U^{2j}E_n^q)''
\to 0$$
where the left-hand term is defined to be the
sub-$\mathcal{O}_Y$-module spanned by the basis elements in the top
line above. The images of the basis elements in the bottom line above
form an $\mathcal{O}_Y$-basis of the right-hand term. Moreover, the
map $R - F$ gives rise to a map of short-exact sequences of sheaves of
abelian groups as follows
$$\xymatrix{
{ 0 } \ar[r] &
{ (\operatorname{gr}_U^{2j}E_n^q)' } \ar[r] \ar[d]^{R' - F'} &
{ \operatorname{gr}_U^{2j}E_n^q } \ar[r] \ar[d]^{R - F} &
{ (\operatorname{gr}_U^{2j}E_n^q)'' } \ar[r] \ar[d]^{R'' - F''} &
{ 0 } \cr
{ 0 } \ar[r] &
{ (\operatorname{gr}_U^{2j}E_{n-1}^q)' } \ar[r] &
{ \operatorname{gr}_U^{2j}E_{n-1}^q } \ar[r] &
{ (\operatorname{gr}_U^{2j}E_{n-1}^q)'' } \ar[r] &
{ 0. } \cr
}$$
The map $R'$ is $\mathcal{O}_Y$-linear, annihilates the
basis elements with $s = n - v - 1$, and leaves the remaining basis
elements unchanged, and the map $F'$ is equal to zero. It follows that
$R' - F' = R'$ is surjective, and hence we have the following
short-exact sequence of kernels of the vertical maps in the diagram
above
$$0 \to (\operatorname{gr}_U^{2j}M_n^q)' \to
\operatorname{gr}_U^{2j}M_n^q \to
(\operatorname{gr}_U^{2j}M_n^q)'' \to 0.$$
Moreover, as $n$ varies, the left-hand term is zero as a sheaf of
pro-abelian groups. Indeed, the structure maps are zero. Similarly,
the map $R''$ above is the $\mathcal{O}_Y$-linear map that annihilates
the basis elements with $s = n - v - 1$ and leaves the remaining basis
elements unchanged. The map $F''$ is the $\varphi$-linear map that
annihilates the basis elements with $s = 1$ and that is given on the
remaining basis elements by
$$\begin{aligned}
{} & F\big( dV^s([x_1]^{i_1} \dots [x_r]^{i_r} [\pi]^j d\log x_{m_1}
\dots d\log x_{m_{q-1}}) \big)\cr
{} & = \bar{x}_1^{k_1}\dots\bar{x}_r^{k_r}
dV^{s-1}([x_1]^{i_1'}\dots [x_r]^{i_r'}[\pi]^j
d\log x_{m_1}\dots d\log x_{m_{q-1}}) \cr
\end{aligned}$$
where $i_m=k_mp^{s-1}+i_m'$ with $0\leqslant i_1',\dots,i_r'<p^{s-1}$. It
follows that a local section
$$\omega = \sum_{s,m,i} a_{m,i}^{(s)} \, dV^s([x]_1^{i_1} \dots
[x_r]^{i_r} [\pi]^j d\log x_{m_1} \dots d\log x_{m_{q-1}})$$
of $\operatorname{gr}_U^{2j}E_n''{}^{q}$ lies in
$\operatorname{gr}_U^{2j}M_n''{}^q$ if and only if the local sections
$a_{m,i}^{(s)}$ of $\mathcal{O}_Y$ satisfy the following system of
equations
$$a_{m,i}^{(s - 1)} = \sum_{0 \leqslant k_1,\dots,k_r < p}
(a_{m,kp^{s-1}+i}^{(s)})^p\bar{x}_1^{k_1}\dots \bar{x}_r^{k_r}.$$
Here $1 < s < n - v$ and the multi-indices $m$ and $i$ vary as in the 
statement of Thm.~\ref{Ubasis}~(ii). We note that a solution to this
system of equations determines and is determined by the local sections
defined by the formula
$$a_m = \sum_{i} (a_{m,i}^{(s)})^{p^s}
\bar{x}_1^{i_1} \dots \bar{x}_r^{i_r}.$$
In this sum, the multi-index $i$ ranges as in the statement of
Thm.~\ref{Ubasis}~(ii) depending on the multi-index $m$, which is
fixed, and the index $1 \leqslant s < n - v$, which is arbitrary but
fixed. It follows that the restriction map
$$R'' \colon (\operatorname{gr}_U^{2j}M_{n+1}^q)'' \hskip5pt \to
(\operatorname{gr}_U^{2j}M_n^q)''$$
is an isomorphism, if $n > v + 1$, and that
$(\operatorname{gr}_U^{2j}M_n^q)''$ is zero, if $n \leqslant v +
1$. Hence, the statement~(ii) is equivalent to the statement that, for
$n > v + 1$, the maps
$$\bar{\rho}_{2j} \colon \Omega_Y^{q-1}/B\Omega_Y^{q-1} \to
(\operatorname{gr}_U^{2j}M_n^q)'' \hskip6mm
(\text{resp.}\ \bar{\rho}_{2j} \colon \Omega_Y^{q-1}/Z\Omega_Y^{q-1}
\to (\operatorname{gr}_U^{2j}M_n^q)'')$$
induced by $\rho_{2j}$ are isomorphisms of sheaves of abelian groups.
To prove this, we note that, for every positive integer $s$,
$\Omega_Y^{q-1}/B\Omega_Y^{q-1}$
(resp.~$\Omega_Y^{q-1}/Z\Omega_Y^{q-1}$) has a canonical structure
of a locally free $\smash{\mathcal{O}_Y^{p^s}}$-module. If
$x_1,\dots,x_r$ are the local coordinates of a neighborhood of $X$
that we considered above, then, in the corresponding neighborhood of
$Y$, an $\smash{\mathcal{O}_Y^{p^s}}$-basis is given by the local
sections
$$\bar{x}_1^{i_1}\dots\bar{x}_r^{i_r}
d\log\bar{x}_{m_1}\dots d\log\bar{x}_{m_{q-1}}$$
where $0 \leqslant i_1, \dots, i_r<p^s$ (resp.~$0 \leqslant i_1,
\dots, i_r < p^s$, not all divisible by $p$), and where $1 \leqslant
m_1 < \dots < m_{q-1} \leqslant r$ are such that, if $m$ is largest
with $i_m$ prime to $p$, then $m_i \neq m$, for all $1 \leqslant i
\leqslant q-1$. We note that, for $1 \leqslant s < n - v$ fixed, the
multi-indices $i$ and $m$ vary in the same way as in the statement of
Thm.~\ref{Ubasis}~(ii). Let $\omega$ be a local section of the sheaf
$\Omega_Y^{q-1}/B\Omega_Y^{q-1}$
(resp.~$\Omega_Y^{q-1}/Z\Omega_Y^{q-1}$). Then, for every positive
integer $s$, we can write $\omega$ as an linear combination
$$\omega = \sum_{m,i} (a_{m,i}^{(s)})^{p^s} \bar{x}_1^{i_1} \dots
\bar{x}_r^{i_r} d\log\bar{x}_{m_1} \dots d\log\bar{x}_{m_{q-1}}$$
with respect to this basis. Then the coefficients $a_{m,i}^{(s)}$
are local sections of $\mathcal{O}_Y$ and, as the index $1 \leqslant s
< n - v$ varies, constitute a solution to the system of equations that
define the subsheaf $\operatorname{gr}_U^{2j}M_n''{}^q \subset
\operatorname{gr}_U^{2j}E_n^q$.
Moreover, Lemma~\ref{dlog} shows that
$$\bar{\rho}_{2j}(\omega) = \sum_{s,m,i} a_{m,i}^{(s)}
dV^s([x_1]^{i_1} \dots [x_r]^{i_r} [\pi]^j d\log x_{m_1} \dots d\log
x_{m_{q-1}})$$
and hence $\bar{\rho}_{2j}$ is an isomorphism as stated. The proof of
the statement~(iii) is completely analogous, and statement~(iv)
follows from Thm.~\ref{Ubasis}~(iv).
\end{proof}

\begin{addendum}\label{coh}The canonical projection
$$i^*(\bar{W}_{\boldsymbol{\cdot}}\Omega_{(X,M_X)}^q)_{F=1} \to
(\bar{W}_{\boldsymbol{\cdot}}\Omega_{(Y,M_Y)}^q)_{F=1}$$
is an isomorphism of pre-sheaves of pro-abelian groups on the small
\'{e}tale site of $Y$. Moreover, the associated sheaf for the
\'{e}tale topology is zero.
\end{addendum}

\begin{proof}We recall from Lemma~\ref{R-F} that the map
$$\operatorname{Fil}_U^2i^*\bar{W}_{\boldsymbol{\cdot}}\Omega_{(X,M_X)}^q
\xrightarrow{1-F}
\operatorname{Fil}_U^2i^*\bar{W}_{\boldsymbol{\cdot}}\Omega_{(X,M_X)}^q$$
is surjective. The isomorphism of the statement now follows from
Lemma~\ref{reduction}. Finally, by the proof of
Addendum~\ref{drwlocalization}, there is split-exact sequence
$$0 \to \bar{W}_{\boldsymbol{\cdot}}\Omega_Y^q \to
\bar{W}_{\boldsymbol{\cdot}}\Omega_{(Y,M_Y)}^q \to
\bar{W}_{\boldsymbol{\cdot}}\Omega_Y^{q-1} \to 0,$$
and we have from~\cite[Prop.~I.3.26]{illusie} that, for the \'{e}tale
topology, the map $1-F$ induces surjections of the left and right-hand
terms.
\end{proof}

\begin{theorem}\label{mainthm}There is a natural exact sequence
$$0 \to i^*R^qj_*\mu_p^{\otimes q} \to
i^*\bar{W}_{\boldsymbol{\cdot}}\Omega_{(X,M_X)}^q \xrightarrow{1-F}
i^*\bar{W}_{\boldsymbol{\cdot}}\Omega_{(X,M_X)}^q \to 0$$
of sheaves of pro-abelian groups on the small \'{e}tale site of $Y$ in
the \'{e}tale topology.
\end{theorem}

\begin{proof}We follow Bloch and Kato~\cite{blochkato} and construct
the left-hand map of the statement by means of the symbol maps
$$i^*R^qj_*\mu_p^{\otimes q} \leftarrow
i^*(M_X^{\operatorname{gp}})^{\otimes q}
\to i^*(\bar{W}_{\boldsymbol{\cdot}}\Omega_{(X,M)}^q)^{F=1}$$
The right-hand map takes a local section $a_1 \otimes \dots \otimes
a_q$ to $d\log a_1\dots d\log a_q$ and the left-hand map takes the
same local section to the symbol $\{a_1,\dots,a_q\}$. We recall the
definition of the latter. By Hilbert's Theorem 90, the Kummer sequence
$$0 \to \mu_p \to \mathcal{O}_{U}^* \xrightarrow{p} \mathcal{O}_{U}^*
\to 0$$
gives rise to an exact sequence
$$0 \to i^*j_*\mu_p \to i^*j_*\mathcal{O}_{U}^* \xrightarrow{p}
i^*j_*\mathcal{O}_{U}^* \xrightarrow{\partial} i^*R^1j_*\mu_p \to 0$$
of sheaves of abelian groups on the small \'{e}tale site of $Y$ in the
\'{e}tale topology. The symbol $\{a\}$ is defined as the image of the
local section $a$ by the composite
$$i^*M_X^{\operatorname{gp}} \xrightarrow{\sim}
i^*j_*\mathcal{O}_{U}^* \xrightarrow{\partial}
i^*R^1j_*\mu_p$$
and $\{a_1,\dots,a_q\}$ as the product of $\{a_1\},\dots,\{a_q\}$.

We may assume that the scheme $X$ is connected. Let $\mathcal{V}'$ be
the strictly henselian local ring of $X$ at the generic point of $Y$,
and let $\mathcal{K}'$ be the quotient field of $\mathcal{V}'$. We let
$\tau \colon \operatorname{Spec}\kappa \to Y$ be the inclusion of the
generic point and consider the maps induced by the symbol maps
$$\tau^*i^*R^qj_*\mu_p^{\otimes q} \leftarrow
\tau^*i^*(M_X^{\operatorname{gp}})^{\otimes q}/p \to
\tau^*(i^*\bar{W}_{\boldsymbol{\cdot}}\Omega_{(X,M_X)}^q)^{F=1}.$$
The left and right-hand terms are canonically isomorphic to the
skyscraper sheafs associated with the pro-abelian groups
$\bar{K}_q^M(\mathcal{K}')$ and $(\bar{W}_{\boldsymbol{\cdot}}
\Omega_{(\mathcal{V}',M_{\mathcal{V}'})}^q)^{F=1}$, respectively. It
follows that the left-hand map is a surjection whose kernel is equal
to the subsheaf generated by the sections $a_1 \otimes \dots \otimes
a_q$ with some $a_i+a_j=1$. By Prop.~\ref{steinberg}, these sections
are annihilated by the right-hand symbol map, so we have an induced
map
\begin{equation}\label{inducedmap}
\tau^*i^*R^qj_*\mu_p^{\otimes q} \to 
\tau^*(i^*\bar{W}_{\boldsymbol{\cdot}}\Omega_{(X,M_X)}^q)^{F=1}.
\end{equation}
This map preserves $U$-filtrations and~\cite[Cor.~1.4.1]{blochkato}
and Thm.~\ref{gr} show that the induced map of filtration quotients is
an isomorphism. It follows that the map is an isomorphism.

We consider the following commutative diagram.
$$\xymatrix{
{i^*R^qj_*\mu_p^{\otimes q}} \ar[d] &
{i^*(M_X^{\operatorname{gp}})^{\otimes q}} \ar[l] \ar[r]  &
{(i^*\bar{W}_{\boldsymbol{\cdot}}\Omega_{(X,M_X)}^q)^{F=1}} \ar[d] \cr
{\tau_*\tau^*i^*R^qj_*\mu_p^{\otimes q}\;} 
\ar[rr]_(.45){\sim}^(.45){~(\ref{inducedmap})} &
{} &
{\tau_*\tau^*(i^*\bar{W}_{\boldsymbol{\cdot}}\Omega_{(X,M_X)}^q)^{F=1}.} \cr
}$$
It is proved in~\cite[Prop.~6.1(i)]{blochkato} that the left-hand
vertical map is injective and in~\emph{op.~cit.}, Cor.~6.1.1, that
the upper left-hand horizontal map is surjective. Moreover, the
right-hand vertical map is injective, since, Zariski locally on $Y$,
the sheaf $i^*\bar{W}_{\boldsymbol{\cdot}}\Omega_{(X,M_X)}^q$ is a
quasi-coherent $\mathcal{O}_Y$-module. It follows that the upper
horizontal maps have the same kernel, and hence the symbol maps give
rise to a map
$$i^*R^qj_*\mu_p^{\otimes q} \to
(i^*\bar{W}_{\boldsymbol{\cdot}}\Omega_{(X,M_X)}^q)^{F=1}.$$
Again this map preserves $U$-filtration,
and~\cite[Cor.~1.4.1]{blochkato} (see also~\cite[Prop.~2.4.1]{tsuji})
and Thm.~\ref{gr} show that the induced map of the associated graded
sheaves is an isomorphism. It follows that the map is an isomorphism
as stated.
\end{proof}

\begin{remark}It is possible from the proof of Thm.~\ref{gr} to
derive the following more precise statement about the injectivity of
the map
$$i^*R^qj_*\mu_p^{\otimes q} \to i^*\bar{W}_n \Omega_{(X,M_X)}^q.$$
As in the statement of Thm.~\ref{Ubasis}, let $v=v(j)$ be the
unique integer such that
$$e\big(\frac{p^{-v}-1}{p^{-1}-1}\big) \leqslant j <
e\big(\frac{p^{-(v+1)}-1}{p^{-1}-1}\big).$$
Then the map is injective, if $n - 1 > v$ for all $0 \leqslant j < e'$.
This, in turn, holds if and only if $p^{n-1} > e'$.
\end{remark}

\begin{proof}[Proof of Thm.~\ref{main}]The surjectivity of $1-F$ is
an immediate consequence of Addendum~\ref{coh}. We show by induction
on $v\geqslant 1$ that the symbol maps
$$i^*R^qj_*\mu_{p^v}^{\otimes q} \leftarrow
i^*(M_X^{\operatorname{gp}})^{\otimes q}/p^v \to
(i^*W_{\boldsymbol{\cdot}}\Omega_{(X,M_X)}^q/p^v)^{F=1}$$
are surjective and have the same kernel. The case $v=1$ is
Thm.~\ref{mainthm}. In the induction step we consider the following
diagram with exact rows, where we have abbreviated
$E_{\boldsymbol{\cdot}}^q=i^*W_{\boldsymbol{\cdot}} \Omega_{(X,M_X)}^q$.
$$\xymatrix{
{} &
{i^*R^qj_*\mu_p^{\otimes q}} \ar[r] &
{i^*R^qj_*\mu_{p^v}^{\otimes q}} \ar[r] &
{i^*R^qj_*\mu_{p^{v-1}}^{\otimes q}} &
{} \cr
{} &
{i^*(M_X^{\operatorname{gp}})^{\otimes q}/p} \ar[r] \ar[u] \ar[d] &
{i^*(M_X^{\operatorname{gp}})^{\otimes q}/p^v} \ar[r] \ar[u] \ar[d] &
{i^*(M_X^{\operatorname{gp}})^{\otimes q}/p^{v-1}} \ar[r] \ar[u] \ar[d] &
{0} \cr
{0} \ar[r] &
{(E_{\boldsymbol{\cdot}}^q/p)^{F=1}} \ar[r] &
{(E_{\boldsymbol{\cdot}}^q/p^v)^{F=1}} \ar[r] &
{(E_{\boldsymbol{\cdot}}^q/p^{v-1})^{F=1}} \ar[r] &
{0} \cr
}$$
The exactness of the lower row follows from Prop.~\ref{highertorsion}
and Addendum~\ref{coh}. By induction the right-hand vertical maps are
surjective and have the same kernel. The same is true for the
left-hand vertical maps. It follows that the middle vertical maps and
the upper right-hand horizontal map are surjective. We claim that the
upper left-hand horizontal map is injective. Indeed, this is
equivalent, by the long-exact cohomology sequence, to the statement
that in the sequence
$$i^*R^{q-1}j_*\mu_p^{\otimes q} \to 
i^*R^{q-1}j_*\mu_{p^v}^{\otimes q} \to
i^*R^{q-1}j_*\mu_{p^{v-1}}^{\otimes q}$$
the right-hand map is surjective. But the cup product by a primitive
$p^v$th root of unity defines an isomorphism of the sheaves
$\mu_{p^v}^{\otimes (q-1)}$ and $\mu_{p^v}^{\otimes q}$ and we have
already proved that the following sequence is exact
$$i^*R^{q-1}j_*\mu_p^{\otimes (q-1)} \to 
i^*R^{q-1}j_*\mu_{p^v}^{\otimes (q-1)} \to
i^*R^{q-1}j_*\mu_{p^{v-1}}^{\otimes (q-1)} \to 0.$$
It remains to show that the middle vertical maps in the diagram above
have the same kernel. To this end, we assume, as in the proof of
Thm.~\ref{mainthm}, that $X$ is connected and let $\tau \colon
\operatorname{Spec} k \to Y$ be the inclusion of the generic point. We
consider the symbol maps
$$\tau^*i^*R^qj_*\mu_{p^v}^{\otimes q} \leftarrow
\tau^*i^*(M_X^{\operatorname{gp}})^{\otimes q}/p^v \to
\tau^*(i^*W_{\boldsymbol{\cdot}}\Omega_{(X,M_X)}^q/p^v)^{F=1}.$$
The kernel of the left-hand map is generated by the symbols
$\{a_1,\dots,a_q\}$ with some $a_i+a_j=1$, and Prop.~\ref{steinberg}
shows that these sections are contained in the kernel of the
right-hand map. Hence, we have an induced map
$$\tau^*i^*R^qj_*\mu_{p^v}^{\otimes q} \to
\tau^*(i^*W_{\boldsymbol{\cdot}}\Omega_{(X,M_X)}^q/p^v)^{F=1}.$$
It is an isomorphism by induction and by the fact that the upper and
lower horizontal rows in the diagram above are short-exact. We
consider the following diagram with exact rows.
$$\xymatrix{
{0} \ar[r] &
{i^*R^qj_*\mu_p^{\otimes q}} \ar[r] \ar[d] &
{i^*R^qj_*\mu_{p^v}^{\otimes q}} \ar[r] \ar[d] &
{i^*R^qj_*\mu_{p^{v-1}}^{\otimes q}} \ar[r] \ar[d] &
{0} \cr
{0} \ar[r] &
{\tau_*\tau^*i^*R^qj_*\mu_p^{\otimes q}} \ar[r] &
{\tau_*\tau^*i^*R^qj_*\mu_{p^v}^{\otimes q}} \ar[r] &
{\tau_*\tau^*i^*R^qj_*\mu_{p^{v-1}}^{\otimes q}.} &
{} \cr
}$$
By induction and by~\cite[Prop.~6.1(i)]{blochkato}, the right
and left-hand vertical maps are injective. Hence, also the middle
vertical map is injective. A similar argument shows that also the
right-hand vertical map in the following diagram is injective.
$$\xymatrix{
{i^*R^qj_*\mu_{p^v}^{\otimes q}} \ar@{^{(}->}[d] &
{i^*(M_X^{\operatorname{gp}})^{\otimes q}/p^v} \ar@{->>}[l] \ar@{->>}[r]  &
{(i^*W_{\boldsymbol{\cdot}}\Omega_{(X,M_X)}^q/p^v)^{F=1}} \ar@{^{(}->}[d] \cr
{\tau_*\tau^*i^*R^qj_*\mu_{p^v}^{\otimes q}\;} 
\ar[rr]^(.44){\sim} &
{} &
{\tau_*\tau^*(i^*W_{\boldsymbol{\cdot}}\Omega_{(X,M_X)}^q/p^v)^{F=1}.} \cr
}$$
It follows that the upper horizontal maps have the same kernel as
desired. This way we obtain the left-hand map of the statement of
Thm.~\ref{main}
$$i^*R^qj_* \mu_{p^v}^{\otimes q} \to
(i^*W_{\boldsymbol{\cdot}} \Omega_{(X,M_X)}^q/p^v)^{F = 1}.$$
Finally, an induction argument based on the short-exactness of the
upper and lower horizontal sequence in the diagram at the beginning of
the proof shows that this map is an isomorphism.
\end{proof}

\section{Henselian discrete valuation rings}\label{bl}

\subsection{}In this section we prove Thm.~\ref{ktheory} of the
introduction. The proof uses the following commutative diagram of
pro-abelian groups in which the right-hand vertical map is the
cyclotomic trace of~\cite{bhm}. We refer the reader to~\cite{h2} for
an introduction and a comprehensive list of references to this
construction.
\begin{equation}\label{diagram}
\xymatrix{
{K_*^M(\mathcal{K})\otimes S_{\mathbb{Z}/p^v}(\mu_{p^v})} \ar[r] \ar[d] &
{K_*(\mathcal{K},\mathbb{Z}/p^v)} \ar[d] \cr
{W_{\boldsymbol{\cdot}}\Omega_{(\mathcal{V},M_{\mathcal{V}})}^*\otimes
  S_{\mathbb{Z}/p^v}(\mu_{p^v})} \ar[r] &
{\operatorname{TR}_*^{\,\boldsymbol{\cdot}}(\mathcal{V}|
\mathcal{K};p,\mathbb{Z}/p^v).} \cr 
}
\end{equation}
We recall from~\cite[Thm.~C]{hm2} and~\cite[Thm.~E]{hm3} that in this
diagram the lower horizontal map is an isomorphism of pro-abelian
groups.

Suppose first that the residue field $\kappa$ is separably
closed. Then Thm.~\ref{main} shows that the left-hand vertical map is
injective and an isomorphism onto the Frobenius fixed set of the
target. Similarly, we show in Prop.~\ref{ktc} below that the the
right-hand vertical map is injective and an isomorphism onto the
Frobenius fixed set of the target. This proves Thm.~\ref{ktheory} in
this case.

In the general case the vertical maps in~(\ref{diagram}) are not
injective, but they still induce surjections onto the Frobenius fixed
sets of the respective targets. Hence, to prove Thm.~\ref{ktheory}, 
we must show that the upper horizontal map induces an isomorphism of
the kernel of the left-hand vertical map onto the kernel of the
right-hand vertical map. We first express the two kernels in terms of
de~Rham-Witt groups and then show that the map in question is an
isomorphism. The proof of the latter occupies most of the section.

\subsection{}Let the field $\mathcal{K}$ be as in the statement of
Thm.~\ref{ktheory}. We first consider the left-hand vertical map
in~(\ref{diagram}).

\begin{prop}\label{milnorktheory}Suppose that $\mu_p\subset K$. Then
there is a natural exact sequence of pro-abelian groups
$$0 \to (\bar{W}_{\boldsymbol{\cdot}}
\Omega_{(\mathcal{V},M_{\mathcal{V}})}^{q-1}\otimes\mu_p)_{F=1} \to
\bar{K}_q^M(\mathcal{K}) \to
(\bar{W}_{\boldsymbol{\cdot}}
\Omega_{(\mathcal{V},M_{\mathcal{V}})}^q)^{F=1} \to 0,$$ 
where the left-hand map takes the class of\, $[a]d\log x_1\dots d\log
x_{q-1}\otimes\zeta$ to the class of the symbol
$\{1+a(1-\zeta)^p,x_1,\dots,x_{q-1}\}$.
\end{prop}

\begin{proof}It follows from~\cite[Thm.~2(1)]{kato1} and from
Thm.~\ref{gr} above that the map that to $\{a_1,\dots,a_q\}$
associates $d\log a_1\dots d\log a_q$ induces an isomorphism of
pro-abelian groups
$$\bar{K}_q^M(\mathcal{K}) /
\operatorname{Fil}_U^{2e'}\bar{K}_q(\mathcal{K}) \xrightarrow{\sim}
(\bar{W}_{\boldsymbol{\cdot}}
\Omega_{(\mathcal{V},M_{\mathcal{V}})}^q)^{F=1}.$$
Indeed, the right-hand side is the stalk at the generic point of $Y$
of the sheaf of pro-abelian groups
$(i^*\bar{W}_{\boldsymbol{\cdot}}\Omega_{(X,M_X)}^q)^{F=1}$ on the
small \'{e}tale site of $Y$ in the Nisnevich
topology. Similarly,~\cite[Thm.~2(1)]{kato1} and Addendum~\ref{coh}
shows that the left-hand map of the statement induces an isomorphism
of pro-abelian groups
$$(\bar{W}_{\boldsymbol{\cdot}}
\Omega_{(\mathcal{V},M_{\mathcal{V}})}^{q-1}\otimes\mu_p)_{F=1}
\xrightarrow{\sim}
\operatorname{Fil}_U^{2e'}\bar{K}_q^M(\mathcal{K}).$$
This completes the proof.
\end{proof}

\begin{remark}Suppose that $\mu_{p^v}\subset K$. One can deduce from
Prop.~\ref{milnorktheory} that there exists a natural exact
sequence of pro-abelian groups
$$0 \to (W_{\boldsymbol{\cdot}}\Omega_{(\mathcal{V},M_{\mathcal{V}})}^{q-1}\otimes\mu_{p^v})_{F=1} \to
K_q^M(\mathcal{K})/p^v \to
(W_{\boldsymbol{\cdot}}\Omega_{(\mathcal{V},M_{\mathcal{V}})}^q/p^v)^{F=1} \to 0.$$
However, we do not have a purely algebraic proof of this deduction. We
also do not have an explicit description of the left-hand map for $v >
1$.
\end{remark}

We now turn our attention to the right-hand vertical map
in~(\ref{diagram}). To this end, we consider the cyclotomic trace map
$$\operatorname{tr} \colon K_q(\mathcal{K},\mathbb{Z}/p) \to
\operatorname{TC}_q^{\,\boldsymbol{\cdot}}
(\mathcal{V}|\mathcal{K};p,\mathbb{Z}/p)$$
from $K$-theory to topological cyclic homology;
see~\cite[Sect.~1]{hm2}. The right-hand side is related to
$\operatorname{TR}_*^{\,\boldsymbol{\cdot}}
(\mathcal{V}|\mathcal{K};p,\mathbb{Z}/p)$ by a natural exact sequence
of pro-abelian groups
$$0 \to \operatorname{TR}_{q+1}^{\,\boldsymbol{\cdot}}
(\mathcal{V}|\mathcal{K};p,\mathbb{Z}/p)_{F=1} \!\xrightarrow{\delta}
\operatorname{TC}_q^{\,\boldsymbol{\cdot}}
(\mathcal{V}|\mathcal{K};p,\mathbb{Z}/p) \to
\operatorname{TR}_q^{\,\boldsymbol{\cdot}}
(\mathcal{V}|\mathcal{K};p,\mathbb{Z}/p)^{F=1}\! \to 0.$$ 
We consider the composition of the left-hand map and the canonical map
$$(\bar{W}_{\boldsymbol{\cdot}}
\Omega_{(\mathcal{V},M_{\mathcal{V}})}^{q+1})_{F=1} \to
\operatorname{TR}_{q+1}^{\,\boldsymbol{\cdot}}
(\mathcal{V}|\mathcal{K};p,\mathbb{Z}/p)_{F=1}.$$

\begin{prop}\label{ktc}For all integers $q$, the cyclotomic trace and
the map $\delta$ give rise to a natural isomorphism of pro-abelian
groups
$$K_q(\mathcal{K},\mathbb{Z}/p) \oplus (\bar{W}_{\boldsymbol{\cdot}}
\Omega_{(\mathcal{V},M_{\mathcal{V}})}^{q+1})_{F=1} \xrightarrow{\sim}
\operatorname{TC}_q^{\,\boldsymbol{\cdot}}
(\mathcal{V}|\mathcal{K};p,\mathbb{Z}/p).$$ 
\end{prop}

\begin{proof}We consider the following diagram of pro-abelian groups,
where the horizontal maps are given by the cyclotomic trace on the
first summand and the boundary map on the second summand, and where
the vertical maps are induced by the canonical projection.
$$\xymatrix{ {\phantom{\operatorname{C}_q^{\,\boldsymbol{\cdot}}}
K_q(\mathcal{V},\mathbb{Z}/p) \oplus (\bar{W}_{\boldsymbol{\cdot}}
\Omega_{\mathcal{V}}^{q+1})_{F=1}} \ar[r] \ar[d] &
{\operatorname{TC}_q^{\,\boldsymbol{\cdot}}
(\mathcal{V};p,\mathbb{Z}/p)} \ar[d] \cr
{\phantom{\operatorname{C}_q^{\,\boldsymbol{\cdot}}}
K_q(\kappa,\mathbb{Z}/p) \oplus (\bar{W}_{\boldsymbol{\cdot}}
\Omega_{\kappa}^{q+1})_{F=1}} \ar[r] &
{\operatorname{TC}_q^{\,\boldsymbol{\cdot}}
(\kappa;p,\mathbb{Z}/p).} \cr }$$
The lower horizontal map is an isomorphism by~\cite[Thm.~4.2.2]{gh},
and we claim that also the top horizontal map is an isomorphism.
Indeed, by Addendum~\ref{coh}, the left-hand vertical map induces an
isomorphism of the second summand of the domain onto the second
summand of the target, so the claim follows from the fact that the map
of relative groups induced by the cyclotomic trace
$$K_q(\mathcal{V},\mathfrak{m}\mathcal{V},\mathbb{Z}/p)
\xrightarrow{\sim} \operatorname{TC}_q^{\,\boldsymbol{\cdot}}
(\mathcal{V},\mathfrak{m}\mathcal{V};p,\mathbb{Z}/p)$$
is an isomorphism of pro-abelian groups. The latter statement, in
turn, is proved in~\cite{mccarthy1},~\cite{suslin1,panin},
and~\cite[Thm.~2.1.1]{gh2}.

We recall that~\cite[Addendum~1.5.7]{hm2} gives a map of localization
sequences
$$\xymatrix{
{\cdots} \ar[r] &
{K_q(\mathcal{V},\mathbb{Z}/p)} \ar[r]^{j_*} \ar[d] &
{K_q(\mathcal{K},\mathbb{Z}/p)} \ar[r]^{\partial} \ar[d] &
{K_{q-1}(\kappa,\mathbb{Z}/p)} \ar[r] \ar[d] &
{\cdots} \cr
{\cdots} \ar[r] &
{\operatorname{TC}_q^{\,\boldsymbol{\cdot}}
(V;p,\mathbb{Z}/p)} \ar[r]^{j_*} &
{\operatorname{TC}_q^{\,\boldsymbol{\cdot}}
(\mathcal{V}|\mathcal{K};p,\mathbb{Z}/p)} \ar[r]^{\partial} &
{\operatorname{TC}_{q-1}^{\,\boldsymbol{\cdot}}
(\kappa;p,\mathbb{Z}/p)} \ar[r] &
{\cdots.} \cr
}$$
Moreover, it follows from Addendum~\ref{drwlocalization} and
Thm.~\ref{gr}(i) that the upper row in the following diagram is
exact. 
$$\xymatrix{
{0} \ar[r] &
{(\bar{W}_{\boldsymbol{\cdot}}
\Omega_{\mathcal{V}}^{q+1})_{F=1}} \ar[r]^{j_*} \ar[d]^{\delta} &
{(\bar{W}_{\boldsymbol{\cdot}}
\Omega_{(\mathcal{V},M_{\mathcal{V}})}^{q+1})_{F=1}}
\ar[r]^{\partial} \ar[d]^{\delta} & 
{(\bar{W}_{\boldsymbol{\cdot}}
\Omega_{\kappa}^q)_{F=1}} \ar[r] \ar[d]^{\delta} &
{0} \cr
{\cdots} \ar[r] &
{\operatorname{TC}_q^{\,\boldsymbol{\cdot}}
(V;p,\mathbb{Z}/p)} \ar[r]^{j_*} &
{\operatorname{TC}_q^{\,\boldsymbol{\cdot}}
(\mathcal{V}|\mathcal{K};p,\mathbb{Z}/p)} \ar[r]^{\partial} &
{\operatorname{TC}_{q-1}^{\,\boldsymbol{\cdot}}
(\kappa;p,\mathbb{Z}/p)} \ar[r] &
{\cdots} \cr
}$$
We claim that the diagram commutes. Indeed, the left-hand square
commutes by the universal property of the de~Rham-Witt complex, and the
description of the upper horizontal map $\partial$ in terms of local
coordinates shows that, in order to show that the right-hand square
commutes, it suffices to show that the lower horizontal map
$\partial$ takes $d\log\pi$ to $1$. But this follows from the
definition of $d\log\pi$ and from the commutativity of the right-hand
square in the previous diagram.

Finally, we combine the two diagrams above as follows. The two
diagrams give rise to a map of long-exact sequences from the sum of
the upper rows in the two diagrams to the common lower row in the two
diagrams. We showed in the beginning of the proof that, in this map of
long-exact sequence, two out of three maps are isomorphisms of
pro-abelian groups. The third map is the map of the statement. This
completes the proof.
\end{proof}

\begin{proof}[Proof of Thm.~\ref{ktheory}]We first note that if the
statement is proved in the basic case $v=1$, then the general case
$v\geqslant 1$ follows inductively by using that the coefficient sequence
breaks up into short-exact sequences
$$0 \to K_q(\mathcal{K},\mathbb{Z}/p) \to
K_q(\mathcal{K},\mathbb{Z}/p^v) \to
K_q(\mathcal{K},\mathbb{Z}/p^{v-1}) \to 0.$$
Hence, it suffices to consider the case $v=1$. It follows from
propositions~\ref{milnorktheory} and~\ref{ktc} that the left and
right-hand vertical maps in~(\ref{diagram}) are surjections onto the
domain and target, respectively, of the canonical map
$$(\bar{W}_{\boldsymbol{\cdot}}\Omega_{(\mathcal{V},M_{\mathcal{V}})}^*
\otimes S_{\mathbb{Z}/p}(\mu_p))^{F=1} \xrightarrow{\sim}
\operatorname{TR}_*^{\,\boldsymbol{\cdot}}
(\mathcal{V}|\mathcal{K};p,\mathbb{Z}/p)^{F=1}.$$
The two propositions further identify the kernel of both of the
vertical maps in~(\ref{diagram}) with the direct sum
$$\bigoplus_{s\geqslant 1}\,(\bar{W}_{\boldsymbol{\cdot}}
\Omega_{(\mathcal{V},M_{\mathcal{V}})}^{q+1-2s} \otimes
\mu_p^{\otimes s})_{F=1}.$$
It remains to show that the map between the two kernels induced by the
upper horizontal map in~(\ref{diagram}) is an isomorphism. This, in
turn, is equivalent to showing that the following
diagram~(\ref{boundarydiagram}) of pro-abelian groups commutes. The
three unmarked maps are as follows: The upper horizontal map is
induced by the lower horizontal map in~(\ref{diagram}); the lower
horizontal map is the composition of the canonical map from Milnor
$K$-theory to algebraic $K$-theory followed by the cyclotomic trace;
and the left-hand vertical map is the left-hand map in
Prop.~\ref{milnorktheory}.
\begin{equation}\label{boundarydiagram}
\xymatrix{
{(\bar{W}_{\boldsymbol{\cdot}}
\Omega_{(\mathcal{V},M_{\mathcal{V}})}^{q-1} \otimes \mu_p)_{F=1}}
\ar[r] \ar[d] &
{\operatorname{TR}_{q+1}^{\,\boldsymbol{\cdot}}
(\mathcal{V}|\mathcal{K};p,\mathbb{Z}/p)_{F=1}} \ar[d]^{\delta} \cr
{\bar{K}_q^M(\mathcal{K})} \ar[r] &
{\operatorname{TC}_q^{\,\boldsymbol{\cdot}}
(\mathcal{V}|\mathcal{K};p,\mathbb{Z}/p) / 
\delta(\bar{W}_{\boldsymbol{\cdot}} 
\Omega_{(\mathcal{V},M_{\mathcal{V}})}^{q+1})_{F=1}} \cr
}
\end{equation}
It follows from Addendum~\ref{coh} that every element of the upper
left-hand term can be written in the form
$[a]d\log x_1\dots d\log x_{q-1}\otimes\zeta$. And since all maps
in~(\ref{boundarydiagram}) are $\bar{K}_*^M(\mathcal{K})$-linear, we
can further assume that $q=1$. Hence, it suffices to prove the
following Prop.~\ref{boundary}.
\end{proof}

\begin{prop}\label{boundary}The image of\, $[a]\otimes\zeta$ by the
composite map
$$W_{\boldsymbol{\cdot}} (\mathcal{V}) \otimes \mu_p \to
\operatorname{TR}_2^{\,\boldsymbol{\cdot}}
(\mathcal{V}|\mathcal{K};p,\mathbb{Z}/p) \xrightarrow{\delta}
\operatorname{TC}_1^{\,\boldsymbol{\cdot}}
(\mathcal{V}|\mathcal{K};p,\mathbb{Z}/p)$$ 
is congruent, modulo
$\delta(\bar{W}_{\boldsymbol{\cdot}}
\Omega_{(\mathcal{V},M_{\mathcal{V}})}^{\,2})_{F=1}$, to
$d\log(1+a(1-\zeta)^p)$.
\end{prop}

\begin{proof}We may assume that the discrete valuation ring
$\mathcal{V}$ is complete. Indeed, the canonical map for $\mathcal{V}$
to the completion of $\mathcal{V}$ induces an isomorphism of all three
terms in the statement. The line of proof is similar to that
of~\cite[Addendum~3.3.9]{hm2}. We apply~\emph{op.~cit.}, Lemma~3.3.10,
to the $3\times 3$-diagram of cofibration sequences
$$\xymatrix{
{E_{11}} \ar[r]^{f_{11}} \ar[d]^{g_{11}} &
{E_{12}} \ar[r]^{f_{12}} \ar[d]^{g_{12}} &
{E_{13}} \ar[r]^{f_{13}} \ar[d]^{g_{13}} &
{\Sigma E_{11}} \ar[d]^{\Sigma g_{11}} \cr
{E_{21}} \ar[r]^{f_{21}} \ar[d]^{g_{21}} &
{E_{22}} \ar[r]^{f_{22}} \ar[d]^{g_{22}} &
{E_{23}} \ar[r]^{f_{23}} \ar[d]^{g_{23}} &
{\Sigma E_{21}} \ar[d]^{\Sigma g_{21}} \cr
{E_{31}} \ar[r]^{f_{31}} \ar[d]^{g_{31}} &
{E_{32}} \ar[r]^{f_{32}} \ar[d]^{g_{32}} &
{E_{33}} \ar[r]^{f_{33}} \ar[d]^{g_{33}} \ar@{}[dr]|-{(-1)} &
{\Sigma E_{31}} \ar[d]^{-\Sigma g_{11}} \cr
{\Sigma E_{11}} \ar[r]^{\Sigma f_{11}} &
{\Sigma E_{12}} \ar[r]^{\Sigma f_{12}} &
{\Sigma E_{13}} \ar[r]^{-\Sigma f_{13}} &
{\Sigma^2 E_{11}} \cr
}$$
obtained as the smash product of the coefficient sequence
$$S^0 \xrightarrow{p} S^0 \to M_p \xrightarrow{\beta} S^1$$
as the first smash factor and the fundamental cofibration sequence
$$\operatorname{TC}^n(\mathcal{V}|\mathcal{K};p) \to
\operatorname{TR}^n(\mathcal{V}|\mathcal{K};p) \xrightarrow{R-F}
\operatorname{TR}^{n-1}(\mathcal{V}|\mathcal{K};p)
\xrightarrow{\delta}
\Sigma\operatorname{TC}^n(\mathcal{V}|\mathcal{K};p)$$
as the second smash factor. We recall from~\emph{op.~cit.},
Lemma~3.3.10, that if we are given classes $e_{ij}\in\pi_*(E_{ij})$
such that $g_{33}(e_{33})=f_{12}(e_{12})$ and
$f_{33}(e_{33})=g_{21}(e_{21})$, then the  sum
$f_{21}(e_{21})+g_{12}(e_{12})$ is in the image of
$\pi_*(E_{11})\to\pi_*(E_{22})$. In the case at hand, we consider the
class
$$e_{33} = [a]_{n-1}^p\cdot b_{\zeta} \in \pi_2(E_{33})=\pi_2(M_p
\wedge \operatorname{TR}^{n-1}(\mathcal{V}|\mathcal{K};p)).$$
We wish to show that the image $e_{31}$ of $e_{33}$ by the map
$$f_{33*} = (\operatorname{id}\wedge\delta)_*  \colon  \pi_2(M_p\wedge
\operatorname{TR}^{n-1}(\mathcal{V}|\mathcal{K};p)) \to
\pi_1(M_p\wedge \operatorname{TC}^n(\mathcal{V}|\mathcal{K};p))$$
is congruent, modulo the image of $\pi_2(E_{23}) \to \pi_1(E_{31})$,
to the class
$$e_{31}' = d\log_n(1+a(1-\zeta)^p).$$
We shall use repeatedly that the canonical map 
$$W_n \Omega_{(\mathcal{V},M_{\mathcal{V}})}^q \to
\operatorname{TR}_q^n(\mathcal{V}|\mathcal{K};p) = 
\pi_q(S^0\wedge\operatorname{TR}^n(\mathcal{V}|\mathcal{K};p))$$
is an isomorphism, if $q\leqslant 2$. This was proved
in~\cite[Thm.~3.3.8]{hm2} for $\mathcal{V}=V$. The general case follows from
this by~\cite[Thms.~B and~C]{hm3}.

By the definition of the Bott element, the image of $e_{33}$ by the map
$$g_{33*}=(\beta\wedge\operatorname{id})_*  \colon 
\pi_2(M_p\wedge\operatorname{TR}^{n-1}(\mathcal{V}|\mathcal{K};p)) \to
\pi_1(S^0\wedge\operatorname{TR}^{n-1}(\mathcal{V}|\mathcal{K};p))$$
is equal to the class
$$e_{13}=[a]_{n-1}^pd\log_{n-1}\zeta.$$
Since we assume that $\mathcal{V}$ is $\mathfrak{m}\mathcal{V}$-adically complete, the
proof of Lemma~\ref{R-F} shows that this class is in the image of 
$$f_{12*}=R-F \colon W_n \Omega_{(\mathcal{V},M_{\mathcal{V}})}^1 \to
W_{n-1} \Omega_{(\mathcal{V},M_{\mathcal{V}})}^1.$$
Indeed, the class $e_{13}$ is contained in
$\operatorname{Fil}_U^{2e''}W_{n-1}
\Omega_{(\mathcal{V},M_{\mathcal{V}})}^1$ and
$e''\geqslant 1$. This also shows that the class $e_{31}$, which we
wish to determine, is contained in the image of the map $g_{21*}
\colon \pi_1(E_{21}) \to \pi_1(E_{31})$. 

We write $\zeta=1+u\pi^{e''}$ with $u\in \mathcal{V}^*$ a unit and
consider the class
$$e_{12} =-\sum_{s=0}^{n-1}\sum_{t=0}^{s}
dV^s([a]_{n-s}^{p^t}[u]_{n-s}^{\phantom{e''}}[\pi]_{n-s}^{e''}).$$

\begin{sublemma}\label{f12e12}$f_{12*}(e_{12}) \in
e_{13} + \operatorname{Fil}_U^{4e''}W_{n-1}
\Omega_{(\mathcal{V},M_{\mathcal{V}})}^1$.
\end{sublemma}

\begin{proof}We assume that $u=1$ (the general case is only
notationally more complicated) and calculate
$$\begin{aligned}
(R-F)e_{12}
{} = & -\sum_{s=0}^{n-2}\sum_{t=0}^{s}
dV^s([a]_{n-1-s}^{p^t}[\pi]_{n-1-s}^{e''}) \cr
{} & +\sum_{s=1}^{n-1}\sum_{t=0}^{s}
dV^{s-1}([a]_{n-s}^{p^t}[\pi]_{n-s}^{e''})
+Fd([a]_n[\pi]_n^{e''}) \cr
{} = & \phantom{+}\sum_{s=0}^{n-2}
dV^s([a]_{n-1-s}^{p^{s+1}}[\pi]_{n-1-s}^{e''}) 
+Fd([a]_n[\pi]_n^{e''}) \cr
{} = & \phantom{+}\sum_{s=0}^{n-2}
d([a]_{n-1}^pV^s([\pi]_{n-1-s}^{e''})) 
+Fd([a]_n[\pi]_n^{e''}) \cr
{} = & \phantom{+}\sum_{s=0}^{n-2}
[a]_{n-1}^p dV^s([\pi]_{n-1-s}^{e''}) \cr
{} & + \sum_{s=0}^{n-2}
p[a]_{n-1}^{p-1}d[a]_{n-1}\cdot V^s([\pi]_{n-1-s}^{e''}) \cr
{} & +[a]_{n-1}^{p-1}[\pi]_{n-1}^{e}d([a]_{n-1}[\pi]_{n-1}^{e''}). \cr
\end{aligned}$$
The summands in the last two lines lie in
$\smash{p\operatorname{Fil}_U^{2e''}=\operatorname{Fil}_U^{2(e+e'')}}$ and the sum in the third
last line is congruent to $[a]_{n-1}^pd\log_{n-1}(1+\pi^{e''})$ modulo
$\operatorname{Fil}_U^{4e''}$ by Lemma~\ref{dlog}. This completes the proof.
\end{proof}

It follows from Lemma~\ref{R-F} that for $m\geqslant 2$,
\begin{equation}\label{add} (R-F)^{-1}(\operatorname{Fil}_U^mW_{n-1}
\Omega_{(\mathcal{V},M_{\mathcal{V}})}^q) =  \operatorname{Fil}_U^mW_n
\Omega_{(\mathcal{V},M_{\mathcal{V}})}^q + \ker(R-F),
\end{equation}
and hence, Sublemma~\ref{f12e12} implies that
$$f_{12*}^{-1}(e_{13} + \operatorname{Fil}_U^{4e''}
W_{n-1}\Omega_{(\mathcal{V},M_{\mathcal{V}})}^1) =
e_{12} + \operatorname{Fil}_U^{4e''}W_n
\Omega_{(\mathcal{V},M_{\mathcal{V}})}^1
+ \operatorname{im}(f_{11*}).$$ 
We next consider the image of this subset by the map
$$g_{12*}=p \colon W_n \Omega_{(\mathcal{V},M_{\mathcal{V}})}^1 \to W_n \Omega_{(\mathcal{V},M_{\mathcal{V}})}^1,$$

\begin{sublemma}The subset $g_{12*}(f_{12*}^{-1}(e_{13} +
\operatorname{Fil}_U^{4e''}W_{n-1}
\Omega_{(\mathcal{V},M_{\mathcal{V}})}^1))$ 
is equal to the subset
$d\log_n(1+a(1-\zeta)^p)+\operatorname{Fil}_U^{2(e'+e'')}W_n
\Omega_{(\mathcal{V},M_{\mathcal{V}})}^1 +
\operatorname{im}(g_{12*}f_{11*}).$
\end{sublemma}

\begin{proof}We again assume $u=1$ and recall from Lemma~\ref{pU} that
$$p\operatorname{Fil}_U^{4e''}W_n
  \Omega_{(\mathcal{V},M_{\mathcal{V}})}^1
  =\operatorname{Fil}_U^{2(e'+e'')}W_n
  \Omega_{(\mathcal{V},M_{\mathcal{V}})}^1.$$
Hence, in view of the equation~(\ref{add}), it will suffice to prove
that
$$g_{12*}(e_{12}) \in d\log_n(1+a(1-\zeta)^p) +
\operatorname{Fil}_U^{2(e'+e'')}W_n
\Omega_{(\mathcal{V},M_{\mathcal{V}})}^1.$$
To this end, we use that in $W_m(\mathcal{V})$,
$$[\pi]_m^{e'} \equiv p(-\sum_{v=0}^{m-1}V^v([\pi]_{m-v}^{e''})),$$
modulo $W_m(\mathfrak{m}^{2e'}\mathcal{V})$. If we rewrite
$$e_{12} =\sum_{s=0}^{n-1}dV^s([a]_{n-s}(-\sum_{v=0}^{n-s-1}
V^v([\pi]_{n-s-v}^{e''}))),$$
this implies that $g_{12*}(e_{12})$ is congruent, modulo
$\operatorname{Fil}_U^{4e'}W_n
\Omega_{(\mathcal{V},M_{\mathcal{V}})}^1$, to the sum
$$\sum_{s=0}dV^s([a]_{n-s}[\pi]_{n-s}^{e'}).$$
Finally, Lemma~\ref{dlog} shows that this sum is congruent, modulo
$\operatorname{Fil}_U^{4e'}W_n
\Omega_{(\mathcal{V},M_{\mathcal{V}})}^1$ to the class
$e_{22} = d\log_n(a(1-\zeta)^p)$.
\end{proof}

Recall that the map
$$\bar{f}_{21*}  \colon  \pi_1(E_{21}) / \operatorname{im}(f_{23*})
\hookrightarrow \pi_1(E_{22})$$
induced by $f_{21*}$ is identified with the canonical inclusion
$$(W_n \Omega_{(\mathcal{V},M_{\mathcal{V}})}^1)^{F=1} \hookrightarrow
W_n \Omega_{(\mathcal{V},M_{\mathcal{V}})}^1.$$
We can now conclude that
$\bar{f}_{21*}^{-1}(g_{12*}(f_{12*}^{-1}(e_{13})))$
is contained in
$$d\log_n(1+a(1-\zeta)^p) +
\operatorname{Fil}_U^{2(e'+e'')}(W_n
\Omega_{(\mathcal{V},M_{\mathcal{V}})}^1)^{F=1}
+ \operatorname{im}(\bar{g}_{11*}).$$
The image of this set by the map
$$\bar{g}_{21*}  \colon  \pi_1(E_{21})/\operatorname{im}(f_{23*}) \to
\pi_1(E_{31})/\operatorname{im}(g_{21*}f_{23*})$$
is equal to the class of $d\log_n(1+a(1-\zeta)^p)$, provided that
$$\operatorname{Fil}_U^{2(e'+e'')}(W_n
\Omega_{(\mathcal{V},M_{\mathcal{V}})}^1)^{F=1} \subset
p(W_n \Omega_{(\mathcal{V},M_{\mathcal{V}})}^1)^{F=1}.$$
We shall prove in Cor.~\ref{multbyp} below that this is almost
true. More precisely, we will show that given $n\geqslant 1$, there exists
$m\geqslant n$ such that the left-hand side is contained in the image of
the composite
$$p(W_m \Omega_{(\mathcal{V},M_{\mathcal{V}})}^1)^{F=1} \;
\hookrightarrow \; W_m \Omega_{(\mathcal{V},M_{\mathcal{V}})}^1
\xrightarrow{R^{m-n}} W_n \Omega_{(\mathcal{V},M_{\mathcal{V}})}^1.$$
The proof of this will occupy the rest of this section. We may then
conclude that given $n\geqslant 1$, there exists $m\geqslant n$ such
that the map of the statement takes $[a]_m\otimes\zeta$ to
$d\log_n(1+a(1-\zeta)^p)$. The proposition follows.
\end{proof}

\begin{lemma}\label{ptorsion}The map that to $x\otimes\zeta$ assigns
$xd\log\zeta$ is an isomorphism of pro-abelian groups
$$W_{\boldsymbol{\cdot}} (\mathfrak{m}\mathcal{V})\otimes\mu_p
\xrightarrow{\sim} (\operatorname{Fil}_U^{2(e''+1)}
W_{\boldsymbol{\cdot}}\Omega_{(\mathcal{V},M_{\mathcal{V}})}^1)[p].$$ 
\end{lemma}

\begin{proof}It follows from~\cite[Thm.~E]{hm3} that the map
$$W_{\boldsymbol{\cdot}} (\mathcal{V})\otimes\mu_p \to
W_{\boldsymbol{\cdot}}\Omega_{(\mathcal{V},M_{\mathcal{V}})}^1[p]$$
that to $x\otimes\zeta$ associates $xd\log\zeta$ is an isomorphism
of pro-abelian groups. This map factors as the composite
$$W_{\boldsymbol{\cdot}} (\mathcal{V})\otimes\mu_p \to
(\operatorname{Fil}_U^{2e''}W_{\boldsymbol{\cdot}}
\Omega_{(\mathcal{V},M_{\mathcal{V}})}^1)[p] \hookrightarrow
W_{\boldsymbol{\cdot}}\Omega_{(\mathcal{V},M_{\mathcal{V}})}^1[p],$$
and since the right-hand map is injective, both maps are necessarily
isomorphisms. We wish to conclude that the map of the
statement is an isomorphism. To this end we consider the following
diagram of pro-abelian group with exact columns.
$$\xymatrix{
{0} \ar[d] &
{0} \ar[d] \cr
{W_{\boldsymbol{\cdot}} (\mathfrak{m}\mathcal{V})\otimes\mu} \ar[r]
\ar[d] &
{(\operatorname{Fil}_U^{2(e''+1)}W_{\boldsymbol{\cdot}}
\Omega_{(\mathcal{V},M_{\mathcal{V}})}^1)[p]} \ar[d] \cr
{W_{\boldsymbol{\cdot}} (\mathcal{V})\otimes\mu} \ar[r]^(.45){\sim}
\ar[d] &
{(\operatorname{Fil}_U^{2e''}W_{\boldsymbol{\cdot}}
\Omega_{(\mathcal{V},M_{\mathcal{V}})}^1)[p]} \ar[d] \cr
{W_{\boldsymbol{\cdot}} (\kappa)\otimes\mu} \ar[r] \ar[d] &
{(\operatorname{Fil}_U^{2e''}W_{\boldsymbol{\cdot}}
\Omega_{(\mathcal{V},M_{\mathcal{V}})}^1 / 
\operatorname{Fil}_U^{2(e''+1)}W_{\boldsymbol{\cdot}}
\Omega_{(\mathcal{V},M_{\mathcal{V}})}^1)[p]} \cr
{0} &
{} \cr
}$$
It suffices to show that the lower horizontal map is
injective. To prove this, we compose with the canonical inclusion
$$(\operatorname{Fil}_U^{2e''}W_{\boldsymbol{\cdot}}
\Omega_{(\mathcal{V},M_{\mathcal{V}})}^1 / 
\operatorname{Fil}_U^{2(e''+1)}W_{\boldsymbol{\cdot}}
\Omega_{(\mathcal{V},M_{\mathcal{V}})}^1)[p]
\hookrightarrow
W_{\boldsymbol{\cdot}}\Omega_{(\mathcal{V},M_{\mathcal{V}})}^1 / 
\operatorname{Fil}_U^{2(e''+1)}W_{\boldsymbol{\cdot}}
\Omega_{(\mathcal{V},M_{\mathcal{V}})}^1$$
and use the isomorphism of Lemma~\ref{reduction}
$$W_{\boldsymbol{\cdot}} \Omega_{(\mathcal{V},M_{\mathcal{V}})}^1 /
\operatorname{Fil}_U^{2(e''+1)}W_{\boldsymbol{\cdot}}
\Omega_{(\mathcal{V},M_{\mathcal{V}})}^1 \xrightarrow{\sim}
W_{\boldsymbol{\cdot}}\Omega_{(\mathcal{V}_{e''+1},M_{e''+1})}^1.$$
The resulting map
$$W_n(\kappa)\otimes\mu_p \to W_{\boldsymbol{\cdot}}
\Omega_{(\mathcal{V}_{e''+1},M_{e''+1})}^1$$
takes $a\otimes\zeta$ to $\tilde{a}d\log_n\bar{\zeta}$, where
$\tilde{a}\in W_n(\mathcal{V}_{e''+1})$ is any lifting of $a\in
W_n(\kappa)$. The ring $\mathcal{V}_{e''+1}$ is isomorphic to the
truncated polynomial ring $\kappa[t]/(t^{e''+1})$, and we can choose
the isomorphism such that the induced map of residue fields is the
identity map. The image of $\bar{\zeta}$ by this isomorphism has the
form $1+ut^{e''}$, where $u\in \kappa^*$ is a unit. (Since $\zeta\in
V$, we can even assume that $u\in k^*$.) Hence, it follows from
Lemma~\ref{dlog} that the composition
$$W_n(\kappa)\otimes\mu_p \to
W_n \Omega_{(\mathcal{V}_{e''+1},M_{e''+1})}^1 \xrightarrow{\sim}
W_n \Omega_{(\kappa[t]/(t^{e''+1}),\mathbb{N}_0)}^1$$
is equal to the $W_n(\kappa)$-linear map that takes $1\otimes\zeta$ to the
sum
$$d\log_n\bar{\zeta}=\sum_{s=0}^{n-1}dV^s([u]_{n-s}[t]_{n-s}^{e''}).$$
The domain of this map is given by
$$W_n(\kappa)/pW_n(\kappa)\otimes\mu_p = W_n(\kappa)/VFW_n(\kappa)
\otimes\mu_p = W_n(\kappa)/VW_n(\kappa^p)\otimes\mu_p,$$
and the target is given by Prop.~\ref{truncated} below. We must
show that if $n$ is sufficiently large, then for all
$a\in W_n(\kappa)$, the product
$$\Theta=ad\log_n\bar{\zeta} = (\sum_{r=0}^{n-1}V^r([a_r]_{n-r}))\cdot
(\sum_{s=0}^{n-1}dV^s([u]_{n-s}[t]_{n-s}^{e''}))$$
is equal to zero if and only if $a_0=0$ and $a_r\in \kappa^p$, for all
$1\leqslant r<n$. We write $e''=p^vi$ with $i$ prime to $p$ and proceed to
rewrite the summands of $\Theta=\sum\Theta_{r,s}$ in the form of 
Prop.~\ref{truncated}. We first note that since $d$ is a derivation
$$\Theta_{r,s} = d(V^r([a_r]_{n-r}) \cdot V^s([u]_{n-s}[t]_{n-s}^{p^vi})) 
-dV^r([a_r]_{n-r}) \cdot V^s([u]_{n-s}[t]_{n-s}^{p^vi}),$$
for all $0\leqslant r,s<n$.

Suppose first that $r>0$. If $0\leqslant s\leqslant v$ and $s\geqslant
r$, we get
$$\begin{aligned}
\Theta_{r,s} 
{} & = p^rdV^s([a_r]_{n-s}^{p^{s-r}}[u]_{n-s}) \cdot [t]_n^{p^{v-s}i} \cr
{} & \hskip3mm
+ p^{r+v-s}iV^s([a_r]_{n-s}^{p^{s-r}}[u]_{n-s}) \cdot 
[t]_n^{p^{v-s}i}d\log_nt \cr
{} & \hskip3mm
- V^s([a_r]_{n-s}^{p^{s-r}}d\log_{n-s}a_r\cdot [u]_{n-s}) \cdot
[t]_n^{p^{v-s}i}. \cr
\end{aligned}$$
The first summand on the right-hand side is zero, since
$$p^rdV^s([a_r]_{n-s}^{p^{s-r}}[u]_{n-s}) \in
\operatorname{Fil}^{s+r}W_n \Omega_{\kappa}^1$$
and $p^{r+s}\cdot p^{s-v}i \geqslant e''+1$, and the second term is zero 
for similar reasons. The third term is zero if and only if
$$V^s([a_r]_{n-s}^{p^{s-r}}d\log_{n-s}a_r\cdot [u]_{n-s}) \in
\operatorname{Fil}^{s+1}W_n \Omega_{\kappa}^1,$$
and this happens if and only if $a_r\in \kappa^p$. Indeed, the filtration
of the groups $W_n\Omega_{\kappa}^q$ is known completely
by~\cite[Prop.~I.2.12]{illusie}. If $0\leqslant s\leqslant v$ and
$s<r$, we have
$$\begin{aligned}
\Theta_{r,s} 
{} & = p^sdV^s(V^{r-s}([a_r]_{n-r})[u]_{n-s}) \cdot 
[t]_n^{p^{v-s}i} \cr
{} & = p^viV^s(V^{r-s}([a_r]_{n-r})[u]_{n-s}) \cdot 
[t]_n^{p^{v-s}i}d\log_nt \cr
{} & \hskip3mm
- V^s(dV^{r-s}([a_r]_{n-r}) [u]_{n-s}) \cdot [t]_n^{p^{v-s}i}, \cr
\end{aligned}$$
and all three terms are zero, since $p^r\cdot p^{v-s}i\geqslant e''+1$. If
$0\leqslant v<s$ and $s\geqslant r$,
$$\begin{aligned}
\Theta_{r,s} 
{} & = p^rdV^{s-v}(V^v([a_r]_{n-s}^{p^{s-r}}[u]_{n-s})
\cdot [t]_{n-s+v}^i) \cr
{} & \hskip3mm
- V^{s-v}(V^v([a_r]_{n-s}^{p^{s-r}}d\log_{n-s}a_r\cdot [u]_{n-s})
\cdot [t]_{n-s+v}^i). \cr
\end{aligned}$$
The first term is zero, since $p^{v+r}i\geqslant e''+1$, and the
second term is equal to zero, if and only if $a_r\in \kappa^p$. Finally, if
$0\leqslant v<s$ and $s<r$, we have
$$\begin{aligned}
\Theta_{r,s} 
{} & = p^sdV^{s-v}(V^v(V^{r-s}([a_r]_{n-r}) [u]_{n-s})
\cdot [t]_{n-s+v}^i) \cr
{} & \hskip3mm
- V^{s-v}(V^v(dV^{r-s}([a_r]_{n-r}) [u]_{n-s})
\cdot [t]_{n-s+v}^i), \cr
\end{aligned}$$
and both terms are zero, since $p^{r-s+v}i\geqslant e''+1$. 

We next evaluate the remaining summands $\Theta_{0,s}$. If $0\leqslant
s\leqslant v$, we have
$$\begin{aligned}
\Theta_{0,s} & = dV^s([a_0]_{n-s}^{p^s}[u]_{n-s}) \cdot
[t]_n^{p^{v-s}i} \cr
{} & \hskip3mm
+p^{v-s}iV^s([a_0]_{n-s}^{p^s}[u]_{n-s})\cdot
[t]_n^{p^{v-s}i}d\log_nt \cr
{} & \hskip3mm
-V^s([a_0]_{n-s}^{p^s}d\log_{n-s}a_0\cdot [u]_{n-s}) \cdot
[t]_n^{p^{v-s}i}. \cr
\end{aligned}$$
The first term is zero if and only if $s=0$ and $a_0\in \kappa^p$, the
second term is zero if and only if $s<v$, and the last term is zero if
and only if $a_0\in \kappa^p$. Finally, if $s>v$,
$$\begin{aligned}
\Theta_{0,s} & = dV^{s-v}(V^v([a_0]_{n-s}^{p^s}[u]_{n-s})\cdot
[t]_{n-s+v}^i) \cr
{} & \hskip3mm -V^{s-v}(V^v([a_0]_{n-s}^{p^s}d\log_{n-s}\cdot
[u]_{n-s}) \cdot [t]_{n-s+v}^i). \cr
\end{aligned}$$
The first term is zero if and only if $a_0$ is zero, and the second
term is zero if and only if $a_0\in \kappa^p$.

We can now show that for $n>v$, the product
$$\Theta = ad\log_n\bar{\zeta}
\in W_n \Omega_{(\kappa[t]/(t^{e''+1}),\mathbb{N}_0)}^1$$
is equal to zero if and only if $a_0=0$ and $a_r\in \kappa^p$, for all
$1\leqslant r<n$, as desired. To this end we use the direct sum
decomposition of the de~Rham-Witt group on the left-hand side
exhibited by Prop.~\ref{truncated} below. Suppose first that $1\leqslant
r<n$. Then $\Theta_{r,s}=0$ if and only if $r>s$ or $r\leqslant s$ and
$a_r\in \kappa^p$. Suppose that $r\leqslant s$. Then the element
$\Theta_{r,s}$ belongs to the direct summand
$V^s(W_{n-s}\Omega_{\kappa}^1) \cdot [t]_n^{p^{v-s}i}$, if $s\leqslant
v$, and to the direct summand $V^{s-v}(W_{n-s+v}\Omega_{\kappa}^q
\cdot [t]_{n-s+v}^i)$, if $s>v$. In particular, two non-zero elements
$\Theta_{r,s}$ and $\Theta_{r',s'}$ belong to the same summand if and
only if $s=s'$. It follows that the sum
$$\sum_{r=1}^{n-1}\sum_{s=0}^{n-1}\Theta_{r,s}
=\sum_{r=1}^{n-1}\sum_{s=r}^{n-1}\Theta_{r,s}$$
is equal to zero if and only if $a_r\in \kappa^p$, for all $1\leqslant
r<n$. A similar argument shows that no cancellation can occur between
the elements $\Theta_{r,s}$, $1\leqslant r\leqslant s<n$, and the
elements $\Theta_{0,s'}$, $0\leqslant s'<n$. Finally, $\Theta_{0,s}$
is non-zero, if $a_0$ is non-zero and $s\geqslant v$. This completes
the proof.
\end{proof}

\begin{cor}\label{multbyp}The map induced from multiplication by $p$
$$p \colon \operatorname{Fil}_U^{2(e''+1)}(W_{\boldsymbol{\cdot}}
\Omega_{(\mathcal{V},M_{\mathcal{V}})}^1)^{F=1} \to
\operatorname{Fil}_U^{2(e'+1)}(W_{\boldsymbol{\cdot}}
\Omega_{(\mathcal{V},M_{\mathcal{V}})}^1)^{F=1}$$
is an isomorphism of pro-abelian groups.
\end{cor}

\begin{proof}We abbreviate
$E_{\boldsymbol{\cdot}}^q = W_{\boldsymbol{\cdot}}
\Omega_{(\mathcal{V},M_{\mathcal{V}})}^q$ and consider the following
diagram.
$$\xymatrix{
{0} \ar[r] &
{\operatorname{Fil}_U^{2(e''+1)}(E_{\boldsymbol{\cdot}}^1)^{F=1}}
\ar[r] \ar[d]^{p} & 
{\operatorname{Fil}_U^{2(e''+1)}E_{\boldsymbol{\cdot}}^1} \ar[r]^{1-F}
\ar[d]^{p} & 
{\operatorname{Fil}_U^{2(e''+1)}E_{\boldsymbol{\cdot}}^1} \ar[r]
\ar[d]^{p} & 
{0} \cr
{0} \ar[r] &
{\operatorname{Fil}_U^{2(e'+1)}(E_{\boldsymbol{\cdot}}^1)^{F=1}}
\ar[r] & 
{\operatorname{Fil}_U^{2(e'+1)}E_{\boldsymbol{\cdot}}^1} \ar[r]^{1-F}
& 
{\operatorname{Fil}_U^{2(e'+1)}E_{\boldsymbol{\cdot}}^1} \ar[r] &
{0} \cr
}$$
It follows from Lemma~\ref{pU} that the middle and right-hand vertical
maps, which are induced by multiplication by $p$, are well-defined and
surjective. The left-hand vertical map is map of kernels induced by
$1-F$. This is the map of the statement. The horizontal maps $1-F$
induce a map between the kernels of the middle and right-hand vertical
maps
$$1-F \colon \operatorname{Fil}_U^{2(e''+1)}
E_{\boldsymbol{\cdot}}^1[p] \to 
\operatorname{Fil}_U^{2(e''+1)}E_{\boldsymbol{\cdot}}^1[p]$$
and Lemma~\ref{ptorsion} shows that this map is an isomorphism of
pro-abelian groups. Indeed, the map
$1-F \colon \bar{W}_{\boldsymbol{\cdot}} (\mathfrak{m}\mathcal{V}) \to
\bar{W}_{\boldsymbol{\cdot}} (\mathfrak{m}\mathcal{V})$ is an
isomorphism, since the geometric series $1+F+F^2+\dots$ converges. 
The corollary follows.
\end{proof}

\appendix

\section{Truncated polynomial algebras}

\subsection{}In this appendix, we give an explicit formula for the
de~Rham-Witt complex of a truncated polynomial algebra in terms of the
de~Rham-Witt complex of the coefficient ring. The formula is derived
from the corresponding formula, proved in~\cite[Thm.~B]{hm3}, for
the de~Rham-Witt complex of a polynomial algebra, and it generalizes
the formula of the thesis of K\aa{}re Nielsen~\cite{nielsen}, where
the case with coefficient ring $\mathbb{F}_p$ was considered.

Let $A$ be a $\mathbb{Z}_{(p)}$-algebra with $p\neq 2$, and let
$A[t]$ be the polynomial algebra in one variable with the pre-log
structure $\alpha \colon \mathbb{N}_0\to A$, $\alpha(i)=t^i$. One can
show as in~\cite[Thm.~B]{hm3} that every element $\omega^{(n)} \in
W_n\Omega_{(A[t],\mathbb{N}_0)}^q$ can be written uniquely
$$\begin{aligned}
\omega^{(n)}  = & \sum_{i\in\mathbb{N}_0}
(a_{0,i}^{(n)}[t]_n^i+b_{0,i}^{(n)}[t]_n^id\log_nt) \cr
{} & + \sum_{s=1}^{n-1} \sum_{i\in I_p} (
V^s(a_{s,i}^{(n-s)}[t]_{n-s}^i) + dV^s(b_{s,i}^{(n-s)}[t]_{n-s}^i) ) \cr
\end{aligned}$$
where $a_{s,i}^{(m)}\in W_m\Omega_A^q$ and $b_{s,i}^{(m)}\in
W_m\Omega_A^{q-1}$, and where $I_p$ denotes the set of positive
integers prime to $p$. The formulas for the product, differential, and
Frobenius and Verschiebung operators may be found in
\emph{op.~cit.}, Sect.~4.2. We now fix an integer $N \geqslant 1$ and
consider the subgroup
$$I_n^q \subset W_n \Omega_{(A[t],\mathbb{N}_0)}^q$$
of those elements $\omega^{(n)}$ such that
$a_{s,i}^{(m)} \in \operatorname{Fil}^vW_m\Omega_A^q$ and
$b_{s,i}^{(m)} \in \operatorname{Fil}^vW_m\Omega_A^{q-1}$, for some
$0 \leqslant v < m$ with $p^vi \geqslant N$. We consider the ring
$A[t]/(t^N)$ with the induced pre-log structure. The following result
expresses $\smash{W_n \Omega_{(A[t]/(t^N),\mathbb{N}_0)}^q}$ as a
direct sum of groups $W_{m-v} \Omega_A^q$ and $W_{m-v} \Omega_A^{q-1}$.

\begin{prop}\label{truncated}The canonical projection induces an
isomorphism
$$W_n \Omega_{(A[t],\mathbb{N}_0)}^q/I_n^q \xrightarrow{\sim}
W_n \Omega_{(A[t]/(t^N),\mathbb{N}_0)}^q.$$
\end{prop}

\begin{proof}We see as in the proof of Lemma~\ref{reduction} that it
suffices to show that $I_n^*$ is a differential graded ideal with
$W_n((t^N)) \subset I_n^0$ and that if $J_n^*$ is another differential
graded ideal with $W_n((t^N)) \subset J_n^0$, then $I_n^*\subset
J_n^*$. We leave the former statement to the reader and prove the
latter. We first show that elements of the form
$V^s(a_{s,i}^{(n-s)}[t]_{n-s}^i)$, where 
$a \in \operatorname{Fil}^vW_{n-s}\Omega_A^q$, for some
$0 \leqslant v < n - s$ with $p^vi \geqslant N$, are contained in
$J_n^q$. By definition of the standard filtration,
$$a_{s,i}^{(n-s)} = V^v(\omega) + dV^v(\omega')$$
for some $\omega \in W_{n-s-v}\Omega_A^q$ and $\omega' \in
W_{n-s-v}\Omega_A^{q-1}$, and hence $V^s(a_{s,i}^{(n-s)} [t]_{n-s}^i)$
is equal to the sum
$$V^{s+v}(\omega[t]_{n-s-v}^{p^vi}) +
p^s dV^{s+v}(\omega'[t]_{n-s-v}^{p^vi}) -
iV^{s+v}(\omega'[t]_{n-s-v}^{p^vi}\,d\log t).$$
We consider the left-hand term. By~\cite[Thm.~A]{hm3}, the
canonical map
$$\Omega_{W_n(A)}^q \to W_n \Omega_A^q$$
is surjective. This shows that $\omega$ can be written as a sum of
elements of the form $x_0 dx_1 \dots dx_q$, where $x_0, \dots, x_q \in
W_{n-s-v}(A)$. But
$$\begin{aligned}
V^{s+v}(x_0 dx_1 \dots dx_q [t]_{n-s-v}^{p^vi})
= V^{s+v}(x_0[t]_{n-s-v}^{p^iv})dV^{s+v}(x_1)\dots dV^{s+v}(x_q),
\end{aligned}$$
which is contained in $J_n^q$, and hence
$V^{s+v}(\omega[t]_{n-s-v}^{p^iv})$, too, is contained in
$J_n^q$. One shows in a similar manner that
$dV^{s+v}(\omega'[t]_{n-s-v}^{p^vi})$ and
$V^{s+v}(\omega'[t]_{n-s-v}^{p^vi}\,d\log t)$ are contained in
$J_n^q$. Hence $V^s(a_{s,i}^{(n-s)}[t]_{n-s}^i)$ is contained in
$J_n^q$. The remaining cases are treated in a completely analogous
manner.
\end{proof}

\section{The Steinberg relation (by Viorel Costeanu)}

\subsection{}This appendix is devoted to the proof of the following
general version of the Steinberg relation in the $p$-typical
de~Rham-Witt complex. We assume as above that the prime $p$ is odd.

\begin{prop}\label{steinberg}Let $(A,M)$ be a
log-$\mathbb{Z}_{(p)}$-algebra, and let $x$ and $y$ be two elements of
$M$ with the property that the sum $\alpha(x) + \alpha(y)$ in $A$ is
equal to $1$. Then the product $d\log_n x \cdot d\log_n y$ in
$W_n\Omega_{(A,M_A)}^2$ is equal to zero.
\end{prop}

\begin{proof}The proof is by induction on $n \geqslant 1$. The proof
of the case $n = 1$ and of the induction step are similar. We assume
that the statement holds for $n - 1$ and prove it for $n$. We have
$$\begin{aligned}
{} & d\log_n x \cdot d\log_n y {} 
 = ([\alpha(x)]_n + [\alpha(y)]_n - 1
- [\alpha(x)]_n - [\alpha(y)]_n) \cdot d\log_nx \cdot d\log_n y \cr
{} & = d[\alpha(x)]_n \cdot d\log_ny + d\log_n x \cdot d[\alpha(y)]_n
+ (1 - [\alpha(x)]_n - [\alpha(y)]_n)d\log_n x \cdot d\log_n y. \cr
\end{aligned}$$
The last summand on the right-hand side is zero by induction. Indeed,
$$1 - [\alpha(x)]_n - [\alpha(y)]_n = [\alpha(x) + \alpha(y)]_n -
[\alpha(x)]_n - [\alpha(y)]_n = V(\tau),$$
for a unique element $\tau \in W_{n-1}(A)$, and hence
$$(1 - [\alpha(x)]_n - [\alpha(y)]_n) \cdot d\log_n x \cdot d\log_n y
= V(\tau \cdot d\log_{n-1} x \cdot d\log_{n-1} y) = 0.$$
It remains to show that $d[\alpha(x)]_n \cdot d\log_n y$ and $d\log_n
x \cdot d[\alpha(y)]_n$ are zero. The proof is the same in the two
cases. Lemma~\ref{congruences} below shows that the polynomial
$$f_s(X) = p^{-s}\big((1 - X)^{p^s} - (1 - X^p)^{p^{s-1}}\big)$$
has integral coefficients. Moreover, one readily verifies the formula
$$[1 - a]_n = [1]_n - [a]_n + \sum_{0 < s < n} V^s(f_s([a]_{n-s}))$$
by evaluating the ghost coordinates of the two sides. It follows that
$$\begin{aligned}
{} & d[\alpha(y)]_n \cdot d\log_n x = d[1-\alpha(x)]_n \cdot d\log_n x \cr 
{} & = dd\log_n x - dd[\alpha(x)]_n +\sum_{0 < s < n}
dV^s(f_s([\alpha(x)]_{n-s})) \cdot d\log_n x. \cr
\end{aligned}$$
The left-hand and middle term in the bottom line are zero. This
completes the proof, for $n = 1$. It remains to prove that the
right-hand term is zero, for $n > 1$. Let $c_{s,i}$ denote the
coefficient of $X^i$ in the polynomial $f_s(X)$. Suppose first that
$p$ does not divide $i$. Since $A$ is assumed to be a
$\mathbb{Z}_{(p)}$-algebra, the integer $i$ is invertible in
$W_n(A)$. Hence, we can write
$$\begin{aligned}
{} & c_{s,i} \cdot dV^s([\alpha(x)]_{n-s}^i) \cdot d\log_n x
 = c_{s,i} \cdot dV^s([\alpha(x)]_{n-s}^i \cdot d\log_{n-s} x) \cr
{} & = i^{-1}c_{s,i} \cdot dV^sd[\alpha(x)]_{n-s} 
 = p^si^{-1}c_{s,i} \cdot ddV^s([\alpha(x)]_{n-s}) \cr
\end{aligned}$$
which is zero as desired. Finally, if $p$ divides $i$, the coefficient
$c_{s,i}$ is divisible by $p^s$ by Lemma~\ref{congruences}. Hence,
we can write
$$\begin{aligned}
{} & c_{s,i} \cdot dV^s([\alpha(x)]_{n-s}^i) \cdot d\log_n x
= p^{-s}c_{s,i} \cdot V^sd([\alpha(x)]_{n-s}^i) \cdot d\log_n x \cr
{} & = p^{-s}ic_{s,i} \cdot V^s([\alpha(x)]_{n-s}^i d\log_{n-s} x)
\cdot d\log_n x \cr
{} & = p^{-s}ic_{s,i} \cdot V^s([\alpha(x)]_{n-s}^i) \cdot d\log_n x
\cdot d\log_n x \cr
\end{aligned}$$
which is zero. This completes the proof.
\end{proof}

\begin{lemma}\label{congruences}The coefficient of $X^i$ in the
polynomial
$$(1 - X)^{p^s} - (1 - X^p)^{p^{s-1}}$$
is divisible by $p^s$, if $p$ does not divide $i$, and by $p^{2s}$, if
$p$ divides $i$.
\end{lemma}

\begin{proof}In general, if two elements $a$ and $b$ of a ring $R$ are
congruent modulo $pR$, then their $p^{s-1}$th powers $a^{p^{s-1}}$ and
$b^{p^{s-1}}$ are congruent modulo $p^sR$. It follows that
$(1-X)^{p^s}$ and $(1-X^p)^{p^{s-1}}$ are congruent modulo
$p^s\mathbb{Z}[X]$. Hence, the coefficients of the polynomial
$(1-X)^{p^s} - (1 - X^p)^{p^{s-1}}$ are divisible by $p^s$. It remains
to show that if $i$ is divisible by $p$, the coefficient of $X^i$ in
this polynomial is divisible by $p^{2s}$. If we write $i = pj$ then
the coefficient in question is equal to
$$\binom{p^s}{pj} - \binom{p^{s-1}}{j} = \binom{p^{s-1}}{j} \big(
(\prod p^s - k)/(\prod k) - 1 \big)$$
where, on the right-hand side, the products range over integers $0 < k
< pj$ that are not divisible by $p$. The $p$-adic valuation of the
first factor on the right-hand side is at least $s - 1 -
v_p(j)$. Indeed, in general, the $p$-adic valuation of the binomial
coefficient $\binom{m+n}{m}$ is equal to the number of carriers in the
addition of the integers $m$ and $n$ in base
$p$~\cite[p.~116]{kummer}. We must show that the $p$-adic valuation of
the second factor on the right-hand side is at least $s + v_p(j) +
1$. Consider the polynomial
$$f(T) = \prod (T - k)$$
where the product ranges over integers $0 < k < pj$ that are not
divisible by $p$. We wish to show that the $p$-adic valuation of
$(f(p^s)/f(0)) - 1$ is at least $s + v_p(j) + 1$, and since $f(0)$ is
not divisible by $p$, this is equivalent to showing that the $p$-adic
valuation of $f(p^s) - f(0)$ is at least $s + v_p(j) + 1$. The
polynomial $f(T) - f(0)$ is divisible by $T$. It suffices to show that
the $p$-adic valuation of the coefficient of $T$ in this polynomial is
at least $v_p(j) + 1$. This coefficient if equal to $f'(0)$. Since
$$d\log f(T) = \sum \frac{dT}{T-k}$$
we find that
$$f'(0) = - f(0) \cdot \sum \frac{1}{k}$$
where the sums range over integers $0 < k < pj$ that are not divisible
by $p$. The number of such integers is equal to $(p-1)j$, and since
$p$ is odd, this is an even number. The partial sum of the $k$th and
$(pj-k)$th summands is equal to
$$\frac{1}{k} + \frac{1}{pj-k} = \frac{pj}{k(pj-k)}$$
which has $p$-adic valuation $v_p(j) + 1$. Finally, the $p$-adic
valuation of a sum of integers is at least the minimum of the $p$-adic
valuations of the summands. Hence the $p$-adic valuation of $f'(0)$ is
at least $v_p(j) + 1$ as desired.
\end{proof}

\providecommand{\bysame}{\leavevmode\hbox to3em{\hrulefill}\thinspace}
\providecommand{\MR}{\relax\ifhmode\unskip\space\fi MR }
\providecommand{\MRhref}[2]{%
  \href{http://www.ams.org/mathscinet-getitem?mr=#1}{#2}
}
\providecommand{\href}[2]{#2}

\end{document}